\documentclass[11pt]{article}
\usepackage{amsthm}
\usepackage{amsmath}
\usepackage{amssymb}
\usepackage[toc,page]{appendix}

\allowdisplaybreaks

\makeatletter
\@addtoreset{equation}{section}

\makeatother

\makeatletter

\newdimen\mainfontsize \mainfontsize=1\@ptsize pt
\topskip=\mainfontsize
   \@maxdepth=\maxdepth

\makeatother

\theoremstyle{plain}
\newtheorem{thm}{Theorem}[section]

\newtheorem{prop}[thm]{Proposition}
\newtheorem{cor}[thm]{Corollary}
\theoremstyle{definition}
\newtheorem{defn}[thm]{Definition}

\theoremstyle{remark}
\newtheorem{rem}[thm]{Remark}

\allowdisplaybreaks
\numberwithin{thm}{section}
\numberwithin{equation}{section}

\title{Fast mean-reversion asymptotics for large portfolios of stochastic volatility models}
\author{Ben Hambly\footnote{hambly@maths.ox.ac.uk} $\,$ and Nikolaos Kolliopoulos\footnote{kolliopoulos@maths.ox.ac.uk (corresponding author)} \\
Mathematical Institute, University of Oxford}
\date{\today} 

\begin{document}

\maketitle

\begin{abstract}
We consider an SPDE description of a large portfolio limit model where the underlying asset prices evolve according to certain stochastic volatility models with default upon hitting a lower barrier. The asset prices and their volatilities are correlated via systemic Brownian motions, and the resulting SPDE is defined on the positive half-space with Dirichlet boundary conditions. We study the convergence of the loss from the system, a function of the total mass of a solution to this stochastic initial-boundary value problem under fast mean reversion of the volatility. We consider two cases. In the first case the volatility converges to a limiting distribution and the convergence of the system is in the sense of weak convergence. On the other hand, when only the mean reversion of the volatility goes to infinity we see a stronger form of convergence of the system to its limit. Our results show that in a fast mean-reverting volatility environment we can accurately estimate the distribution of the loss from a large portfolio by using an approximate constant volatility model which is easier to handle.
\end{abstract}

\section{Introduction}


In this paper our aim is to investigate the fast mean reverting volatility asymptotics for an SPDE-based structural model for portfolio credit. SPDEs arising from large portfolio limits of collections of defaultable constant volatility models were initially studied in Bush, Hambly et al. \cite{BHHJR}, and their regularity was further investigated in Ledger \cite{Ledger14}. In Hambly and Kolliopoulos \cite{HK17,HKE17,HK19} we extended this work to a two-dimensional stochastic volatility setting, and here we consider the question of effective one-dimensional constant volatility approximations which arise by considering fast mean-reversion in the volatilities. This approach is to some extent motivated by the ideas of Fouque, Papanicolaou and Sircar \cite{FOU}, but instead of option prices we look at the systemic risk of large credit portfolios in the fast mean-reverting volatility setting. 

The literature on large portfolio limit models in credit can be divided into two approaches based on either structural or reduced form models for the individual assets. Our focus will be on the structural approach where we assume that we are modelling the financial health of the firms directly and default occurs when these health processes hit a lower barrier. 

The reduced form setting assumes that the default of each firm occurs as a Poisson process and we model the default intensities directly. These can be correlated through systemic factors and through the losses from the portfolio. The evolution of the large portfolio limit of the empirical measure of the loss can be analysed as a law of large numbers and then Gaussian fluctuations derived around this limit, see Giesecke, Sirignano et al. \cite{GSS,SG,SSG,GSS2} and Cvitanic et al. \cite{CMZ}. Further, the large deviations can be analysed, see Sowers and Spiliopoulos \cite{SSO1,SSO2}. It is also possible to take an approach through interacting particle systems where each firm as in one of two states representing financial health and financial distress and there is a movement between states according to some intensity, often firm dependent, and dependent on the proportion of losses, see for instance Dai Pra and Tolotti \cite{PT09} or Dai Pra et al. \cite{PRST}.

Our underlying set up is a structural model for default in which each asset has a distance 
to default, which we think of as the logarithmically scaled asset price process. The asset 
price evolves according to a general stochastic volatility model, in which the distance to 
default of the $i$-th asset
$X^{i}$ satisfies the system
\begin{equation}
\begin{array}{rcl}
dX_{t}^{i} &=& \left(r_{i}-\frac{h^{2}(\sigma_{t}^{i})}{2}\right)dt+h(\sigma_{t}^{i})\left(\sqrt{1-\rho_{1,i}^{2}}dW_{t}^{i}+\rho_{1,i}dW_{t}^{0}\right),\;\;0\leq t\leq T_{i} \\
d\sigma_{t}^{i} &=& k_{i}(\theta_{i}-\sigma_{t}^{i})dt+\xi_{i}g\left(\sigma_{t}^{i}\right)\left(\sqrt{1-\rho_{2,i}^{2}}dB_{t}^{i}+\rho_{2,i}dB_{t}^{0}\right),\;\;t\geq 0 \\
X_{t}^{i} &=& 0,\;t > T_{i} \\
(X_{0}^{i},\,\sigma_{0}^{i}) &=& (x^{i},\sigma^{i,init}),\\
T_i &=& \inf\{ t\geq 0: X_t^{i}=0\}
\end{array}
\label{eq:generalmodel}
\end{equation} 
for all $i \in \mathbb{N}$, where the coefficient vectors $C_i = (r_i, \rho_{1,i}, \rho_{2,i}, k_i, \theta_i, \xi_i)$ are picked randomly and independently from some probability distribution with $\rho_{1,i}, \rho_{2,i} \in [0, \, 1)$, the infinite sequence $\{(x^1, \, \sigma^{1,init}), \, (x^2, \, \sigma^{2,init}), \, ...\}$ of random vectors in $\mathbb{R}^2$ is assumed to be exchangeable, and $g,h$ are functions for which we will give suitable conditions later. The exchangeablity condition implies (see \cite{Aldous,KKP}) the existence of a $\sigma$-algebra $\mathcal{G} \subset \sigma(\{(x^i, \, \sigma^i): \, i \in \mathbb{N}\})$, given which the two-dimensional random vectors $(x^i, \, \sigma^i)$ are pairwise independent and identically distributed. The idiosyncratic Brownian motions $W^i, B^i$ for $i\in \mathbb{N}$ are taken to be pairwise independent, and also independent of the systemic Brownian motions $W^0, B^0$ which have a constant correlation $\rho_3$.

We regard this as a system for $Z^i=(X^i,\sigma^i)$ with 
\[ dZ^i = b^i(Z^i) dt + \Sigma^i(Z^i) d{\mathbf{W}^i}, \;\; Z^i_0 = (x^i,\sigma^{i,init}) \]
for $t<T_i$, where 
\[ b^i(X,\sigma)= (r_i-\frac{h(\sigma)^2}{2}, k_i(\theta_i-\sigma))^{\top},\]
\[ \Sigma^i(X,\sigma) = \left[ 
\begin{array}{cccc}
h(\sigma)\sqrt{1-\rho_{1,i}^{2}} & h(\sigma)\rho_{1,i} & 0 & 0\\
0 & 0 & \xi_{i}g\left(\sigma\right) \sqrt{1-\rho_{2,i}^{2}} & \xi_{i}g\left(\sigma\right) \rho_{2,i} \end{array} \right] \]
and $\mathbf{W}^i = (W^i, W^0, B^i, B^0)^{\top}$. Then, the infinitesimal generator of the above two-dimensional process is given by
\[ \mathcal{A}^{i}f = \sum_{j=1}^2 b_j^{i} \frac{\partial f}{\partial x_j} + \frac12 \sum_{j,k=1}^2 a_{jk}^{i} \frac{\partial^2f}{\partial x_j\partial x_k} \]
for $f\in C^2(\mathbb{R}_{+}\times\mathbb{R}, \mathbb{R})$. The matrix $A^{i} = a_{jk}^{i}$ is given by
\[ A^i = \left[ \begin{array}{cc}
h(\sigma)^2     &  h(\sigma)\xi_ig(\sigma) \rho_{1,i}\rho_{2,i}\rho_3 \\
 h(\sigma)\xi_ig(\sigma) \rho_{1,i}\rho_{2,i}\rho_3     & \xi_i^2 g(\sigma)^2
\end{array} \right], \]
as $A^{i} = \Sigma^{i} R (\Sigma^{i})^{\top}$, with $R$ the covariance matrix for the 4-dimensional Brownian motion $\mathbf{W}^i$.

We can show that the empirical measure of a sequence of finite sub-systems 
\[ \nu^N_t = \frac1{N} \sum_{i=1}^N \delta_{X^i_t,\sigma^i_t}, \]
converges weakly as $N\to\infty$ (see \cite{HK19}) to the probability distribution of $Z_t^1$ given $W^0$, $B^0$ and $\mathcal{G}$. This measure consists of two parts; its restriction to the line $x=0$, which is approximated by the restriction of $\nu^N$ to this line, and its restriction to $\mathbb{R}_{+} \times \mathbb{R}$ which possesses a two-dimensional density $u(t,x,y)$. The density $u(t,x,y)$ can be regarded as an average of solutions to certain two-dimensional SPDEs with a Dirichlet boundary condition on the line $x=0$. In particular, we can write $u = \mathbb{E}[u_{C_1} \, | \, W^0, \, B^0, \, \mathcal{G}]$, where $u_{C_1}(t,x,y)$ is the probability density of $Z_t^1$ given $W^0$, $B^0$, $\mathcal{G}$ and $C_{1}$ on $\mathbb{R}_{+} \times \mathbb{R}$, which satisfies, for any value of the coefficient vector $C_1$, the two-dimensional SPDE
\begin{equation}
du_{C_1} = \mathcal{A}^{1,*} u_{C_1} dt + \mathcal{B}^{1,*}u_{C_1} d(W^0,B^0)^{\top},
\label{maineq}
\end{equation}
where $\mathcal{A}^{1,*}$ is the adjoint of the generator $\mathcal{A}^{1}$ of $Z^1$, and the operator $\mathcal{B}^{1,*}$ is given by
\[ \mathcal{B}^{1,*}f = \left(-\rho_{1,1}h(y)\frac{\partial f}{\partial x},  -\xi_{1}\rho_{2,1}g\left(y\right) \frac{\partial f}{\partial y} \right). \] 
The boundary condition is that $u_{C_1}(t,0,y)=0$ for all $y\in\mathbb{R}$. In the special case where the coefficients are constants independent of $i$, $u$ is itself a solution to the stochastic partial differential equation (\ref{maineq}).

One reason for studying the large portfolio limit is the need to have a useful approximation
which captures the dynamics among the asset prices when the number of assets is large. Moreover, by studying the limit SPDE instead of a finite sub-system of (\ref{eq:generalmodel}), we can potentially provide a more efficient approach to capturing the key drivers of a large portfolio without having to simulate a large number of idiosyncratic Brownian paths.

Of central importance will be the loss function $L$, the mass of the probability distribution of $Z_t^1$ given $W^0$, $B^0$ and $\mathcal{G}$ on the line $x=0$, which measures the total loss in the large portfolio limit. 
The distribution of this function is a simple measure of risk for the portfolio of assets and can be used to find the probability of a large loss, or to determine the prices of portfolio credit derivatives such as CDOs that can be written as expectations of suitable functions of $L$. Thus our focus will be on estimating probabilities of the form
\begin{eqnarray}
\mathbb{P}\Big[L_{t} \in \left(1 - b, \, 1 - a \right)\Big] = \mathbb{P}\Big[\mathbb{P}\big[X_t^{1} > 0 \, \big| \, W^{0}, \, B^{0}, \, \mathcal{G}\big] \in \left(a, \, b \right)\Big]
\label{lossprob}
\end{eqnarray}
for some $0 \leq a < b \leq 1$, that is the probability that the total loss from the portfolio lies within a certain range. Probabilities of the above form can be approximated numerically with a simulated sample of values of $L_t$, obtained via
\begin{eqnarray}\label{averaging}
1 - L_t &=& \mathbb{P}\Big[X_t^{1} > 0 \, \Big| \, W^{0}, \, B^{0}, \, \mathcal{G}\Big] \nonumber \\
&=& \int_{0}^{+\infty}\int_{0}^{+\infty}\mathbb{E}\Big[u_{C_1}(t,x,y) \, \Big| \, W^0, \, B^0, \, \mathcal{G} \Big]dxdy \nonumber \\
&\approx& \frac{1}{n}\sum_{i=1}^{n}{\int_{0}^{+\infty}\int_{0}^{+\infty}u_{c_{1,i}}(t,x,y)dxdy}
\end{eqnarray} 
after solving the SPDE (\ref{maineq}) for $u_{C_1}$ numerically, for a sample $\{c_{1,1}, \, ... , \, c_{1,n} \}$ of values of the vector $C_1$. In the special case when asset prices are modelled as simple constant volatility models, the numerics (see Giles and Reisinger \cite{GR}, or Bujok and Reisinger \cite{BR12} for jump-diffusion models) have a significantly smaller computational cost, which motivates the investigation of the existence of accurate approximations using a constant volatility setting in the general case. We note also that one-dimensional SPDEs describing large portfolio limits in constant volatility environments have been found to have a unique regular solution (see \cite{BHHJR}, or Hambly and Ledger \cite{HL16} for a loss-dependent correlation model), an important component of the numerical analysis and a counterpoint to the fact that we have been unable to establish uniqueness of solutions to the two-dimensional SPDE arising in the CIR volatility case \cite{HK17}. 


We will derive our one-dimensional approximations under two different settings with fast mean-reverting volatility. In what we call the large vol-of-vol setting, 
the mean reversion and volatility in the second equation in (\ref{eq:generalmodel}) are scaled by suitable powers of $\epsilon$ in that
$k_i=\kappa_i/\epsilon$ and $\xi_i = v_i/\sqrt{\epsilon}$ giving
\[ d\sigma_{t}^{i} = \frac{\kappa_{i}}{\epsilon}(\theta_{i}-\sigma_{t}^{i})dt+\frac{v_i}
{\sqrt{\epsilon}} g\left(\sigma_{t}^{i}\right)\left(\sqrt{1-\rho_{2,i}^{2}}dB_{t}^{i}+\rho_{2,i}dB_{t}^{0}\right),\;\;t\geq 0, \]
and then we take $\epsilon \to 0$.
This is distributionally equivalent to speeding up the volatility processes by scaling the time $t$ by $\epsilon$, when $\epsilon$ is small.
Our aim is to take the limit as $\epsilon\to 0$, so that when the system of volatility processes is positive reccurent, averages over finite time intervals involving the sped up volatility processes will approximate the corresponding stationary means. In the limit we obtain a constant volatility large portfolio model which could be used as an effective approximation when volatilities are fast mean-reverting. However, this speeding up does not lead to strong convergence of the volatility processes, allowing only for weak convergence of our system, which can only be established when $\rho_3 = 0$ (effectively separating the time scales) and when $(\kappa_{i},\,\theta_{i},\,v_{i},\,\rho_{2,i})$ is the same constant vector $(\kappa,\,\theta,\,v,\,\rho_{2})$ for all $i \in \mathbb{N}$. 

The case of small vol-of-vol has the mean reversion in the second equation in (\ref{eq:generalmodel}) scaled by $\epsilon$ in that $k_i=\kappa_i/\epsilon$ and
\[ d\sigma_{t}^{i} = \frac{\kappa_{i}}{\epsilon}(\theta_{i}-\sigma_{t}^{i})dt+\xi_{i}
g\left(\sigma_{t}^{i}\right)\left(\sqrt{1-\rho_{2,i}^{2}}dB_{t}^{i}+\rho_{2,i}dB_{t}^{0}\right),\;\;t\geq 0. \]
We regard this case as a small noise perturbation of the constant volatility model, where 
volatilities have stochastic behaviour but are pulled towards their mean as soon as they 
move away from it due to a large mean-reverting drift. When $\epsilon \to 0$, the drifts of 
the volatilities tend to infinity and dominate the corresponding diffusion parts since the 
vol-of-vols remain small, allowing for the whole system to converge to a constant 
volatility setting in a strong sense. 
This strong convergence allows the rate of convergence of probabilities of the form (\ref{lossprob}) to be estimated and gives us a quantitative measure of the loss in accuracy in the estimation of these probabilities when a constant volatility large portfolio model is used to replace a more realistic stochastic volatility perturbation of that model.

In Sections~2 and~3 we present our main results for both settings. The results are then
proved in Sections~4 and~5. Finally the proofs of two propositions showing the positive 
recurrence, and hence applicability of our results, for two classes of models can be 
found in the Appendix.

\section{The main results: large vol-of-vol setting}

We begin with the study of the fast mean-reversion - large vol-of-vol setting, for which we need to assume that the correlation $\rho_3$ of $W^{0}$ and $B^{0}$ is zero. When $g$ is either the square root function or a function behaving almost like a positive constant for large values of the argument, it has been proven in Theorem 4.3 in \cite{HK17} and in Theorem 4.1 in \cite{HK19} respectively that 
\begin{eqnarray*}
&& u_{C_1}(t,\,x,\,y) \\
&& \qquad \qquad = p_{t}\left(y|B^{0},\mathcal{G}\right)
\mathbb{E}\bigg[u\Big(t,x,W^{0},\mathcal{G},C_{1},h\left(\sigma^{1}\right)\Big)\,\bigg|W^{0},\sigma_{t}^{1}= y,B^{0},C_{1},\mathcal{G}\bigg], \nonumber
\end{eqnarray*}
where $p_t$ is the density of each volatility path when the path of $B^{0}$ and the information in $\mathcal{G}$ are given, and $u(t,x,W^{0},\mathcal{G},C_{1},h(\sigma^{1}))$ is the unique $H_0^1\left(0, \, +\infty\right)$ solution to the SPDE 
\begin{eqnarray}\label{1dimeq}
u(t,\,x) &=& u_{0}(x) - \int_{0}^{t}\left(r-\frac{h^2\left(\sigma_{s}^{1}\right)}{2}\right)u_{x}(s,\,x)ds \nonumber \\
&& +\int_{0}^{t}\frac{h^2\left(\sigma_{s}^{1}\right)}{2}
u_{xx}(s,\,x)ds-\rho_{1,1}\int_{0}^{t}h\left(\sigma_{s}^{1}\right)u_{x}(s,\,x)dW_{s}^{0}, 
\end{eqnarray}
where $u_0$ is the density of each $x^i$ given $\mathcal{G}$. In the above expression for the two-dimensional density $u_{C_1}(t,\,x,\,y)$, averaging happens with respect to the idiosyncratic noises, and since we are interested in probabilities concerning $L_t$ which is computed by substituting that density in (\ref{averaging}), averaging happens with respect to the market noise $(W^{0}, \, B^{0})$ as well. Therefore, we can replace $(W^{i}, \, B^{i})$ for all $i \geq 0$ in our system by objects having the same joint law. 
In particular, setting $k_i=\kappa_i/\epsilon$ and $\xi_i = v_i/\sqrt{\epsilon}$, the $i$-th asset's distance to default $X^{i, \epsilon}$ satisfies the system  
\[
\begin{array}{rcl}
X_{t}^{i, \epsilon} &=& x^i + \int_{0}^{t}\left(r_{i}-\frac{h^{2}(\sigma_{t}^{i, \epsilon})}{2}\right)dt \\
&& \qquad \qquad + \int_{0}^{t}h(\sigma_{t}^{i, \epsilon})\left(\sqrt{1-\rho_{1,i}^{2}}dW_{t}^{i}+\rho_{1,i}dW_{t}^{0}\right),\;\;0\leq t\leq T_{i} \\
\sigma_{t}^{i, \epsilon} &=& \sigma^{i,init} + \frac{\kappa_{i}}{\epsilon}\int_{0}^{t}(\theta_{i}-\sigma_{s}^{i, \epsilon})ds + \frac{v_{i}}{\sqrt{\epsilon}}\int_{0}^{t}g\left(\sigma_{s}^{i, \epsilon}\right)d\left(\sqrt{1-\rho_{2,i}^{2}}B_{s}^{i} + \rho_{2,i}B_{s}^{0}\right) \\
X_{t}^{i, \epsilon} &=& 0,\;t > T_{i}^{\epsilon} \\
T_i^{\epsilon} &=& \inf\{ t\geq 0: X_t^{i, \epsilon}=0\}
\end{array} 
\] 
where the $\epsilon$ superscripts are used 
to underline the dependence on $\epsilon$, 
and if we substitute $t = \epsilon t'$ and $s = \epsilon s'$ for $0 \leq s' \leq t'$ and then replace $(W^{i}, \, B^{i})$ by $(W^{i}, \sqrt{\epsilon}B_{\frac{\cdot}{\epsilon}}^{i})$ for all $i \geq 0$ which have the same joint law, the SDE satisfied by the $i$-th volatility process becomes 
\[
\sigma_{\epsilon t'}^{i, \epsilon} = \sigma^{i,init} + \kappa_{i}\int_{0}^{t'}(\theta_{i}-\sigma_{\epsilon s'}^{i, \epsilon})ds' + v_{i}\int_{0}^{t'}g\left(\sigma_{\epsilon s'}^{i, \epsilon}\right)d\left(\sqrt{1-\rho_{2,i}^{2}}B_{s'}^{i} + \rho_{2,i}B_{s'}^{0}\right). 
\] 
This shows that $\sigma^{i,\epsilon} = \sigma_{\epsilon \times \frac{\cdot}{\epsilon}}^{i,\epsilon}$ can be replaced by $\sigma_{\frac{\cdot}{\epsilon}}^{1,1}$ for all $i \geq 1$, i.e the $i$-th volatility process of our model when the mean-reversion coefficient and the vol-of-vol are equal to $\kappa_i$ and $v_i$ respectively and when the time $t$ is scaled by $\epsilon$, speeding up the system of the volatilities when $\epsilon$ is small. 

If $g$ is now chosen so that the system of volatility processes becomes positive recurrent, averages over finite time intervals converge to the corresponding stationary means as the speed tends to infinity, i.e as $\epsilon \to 0^+$, which is the key for the convergence of our system. We give a definition of the required property for $g$.

\begin{defn}[\bf{Positive recurrence property}]\label{reccur}{\rm
We fix the distribution from which each $C_i' = (r_i, \rho_{1,i}, \rho_{2,i}, k_i, \theta_i, \xi_i)$ is chosen, and we denote by $\mathcal{C}$ the $\sigma$-algebra generated by all these coefficient vectors. Then, we say that $g$ has the positive recurrence property when the two-dimensional process $(\sigma_{\cdot}^{i,1}, \, \sigma_{\cdot}^{j,1})$ is a positive recurrent diffusion for any two $i, j \in \mathbb{N}$, for almost all values of $C_i'$ and $C_j'$. This means that given $\mathcal{C}$, there exists a two-dimensional random variable $(\sigma^{i, j, 1, *}, \, \sigma^{i, j, 2, *})$ whose distribution is stationary for $(\sigma_{\cdot}^{i,1}, \, \sigma_{\cdot}^{j,1})$, and whenever $\mathbb{E}[|F(\sigma^{i, j, 1, *}, \, \sigma^{i, j, 2, *})| \, | \, \mathcal{C}]$ exists and is finite for some measurable function $F: \mathbb{R}^2 \rightarrow \mathbb{R}$ we also have:
\[
\lim_{T \rightarrow + \infty}\frac{1}{T}\int_{0}^{T}F\left(\sigma_{s}^{i,1}, \, \sigma_{s}^{j,1}\right)ds = \mathbb{E}\Big[F\left(\sigma^{i, j, 1, *}, \, \sigma^{i, j, 2, *}\right) \, \Big| \, \mathcal{C}\Big],
\]
or equivalently, after a change of variables,
\[
\lim_{\epsilon \rightarrow 0^{+}}\frac{1}{t}\int_{0}^{t}F\left(\sigma_{\frac{s}{\epsilon}}^{i,1}, \, \sigma_{\frac{s}{\epsilon}}^{j,1}\right)ds = \mathbb{E}\Big[F\left(\sigma^{i, j, 1, *}, \, \sigma^{i, j, 2, *}\right) \, \Big| \, \mathcal{C}\Big]
\]
for any $t \geq 0$, $\mathbb{P}$-almost surely.}
\end{defn}

The positive recurrence property is a prerequisite for our convergence results to hold, and now we will state two propositions which give us a few classes of models for which this property is satisfied. The first shows that for the Ornstein-Uhlenbeck model ($g(x) = 1$ for all $x \in \mathbb{R}$) we always have the positive recurrence property. The second shows that for the CIR model ($g(x) = \sqrt{|x|}$ for all $x \in \mathbb{R}$) we have the positive recurrence property provided that the random coefficients of the volatilities satisfy certain conditions. The proofs of both propositions can be found in the Appendix.

\begin{prop}\label{thm3}
Suppose that $g$ is a differentiable function, bounded from below by some $c_g > 0$. Suppose also that $g'(x)\kappa_i(\theta_i - x) < \kappa_ig(x) + \frac{v_i}{2}g''(x)g^2(x)$ for all $x \in \mathbb{R}$ and $i \in \mathbb{N}$, for all possible values of $C_i$. Then $g$ has the positive recurrence property.
\end{prop}

\begin{prop}\label{thm4}
Suppose that $g(x) = \sqrt{|x|}\tilde{g}(x)$, where the function $\tilde{g}$ is a continuously differentiable, strictly positive and increasing function taking values in $[c_g, \, 1]$ for some $c_g > 0$. Then, there exists an $\eta > 0$ such that $g$ has the positive recurrence property when $\Vert C_i - C_j \Vert_{L^{\infty}(\mathbb{R}^6)} < \eta$ and $\frac{\kappa_i}{v_j^2} > \frac{1}{4} + \frac{1}{\sqrt{2}}$ for all $i, \, j \in \mathbb{N}$, $\mathbb{P}$ - almost surely.
\end{prop}

We can proceed now to our main results, which will be governed by the conditional moments $\sigma_{1,1} = \mathbb{E}[h(\sigma^{1,1,1,*}) \, | \, \mathcal{C}]$ and $\sigma_{2,1} = \sqrt{\mathbb{E}[h^2(\sigma^{1,1,1,*}) \, | \, \mathcal{C}]}$, as well as the quantity $\tilde{\sigma} = \sqrt{\mathbb{E}[h(\sigma^{1, 2, 1, *})h(\sigma^{1, 2, 2, *}) \, | \, \mathcal{C}]}$, where $\sigma^{1, 1, 1, *}$, $\sigma^{1, 2, 1, *}$ and $\sigma^{1, 2, 2, *}$ are given in Definition~\ref{reccur}. The next theorem implies the weak convergence of the loss $L_t^{\epsilon} = 1 - \mathbb{P}[X_t^{1, \epsilon} > 0 \, | \, W^{0}, \, B^{0}, \, \mathcal{G}]$ under the fast mean-reverting volatility setting to the loss under an appropriate constant volatility setting.

\begin{thm}\label{theorem}
Suppose that $(\kappa_{i},\,\theta_{i},\,v_{i},\,\rho_{2,i}) = (\kappa,\,\theta,\,v,\,\rho_{2})$ for all $i \in \mathbb{N}$, which is a deterministic vector in $\mathbb{R}^4$, the function $h$ is bounded, and that $g$ has the positive recurrence property, in which case we have $\sigma_{1,1} = \mathbb{E}[h(\sigma^{1,1,1,*})]$, $\sigma_{2,1} = \sqrt{\mathbb{E}[h^2(\sigma^{1,1,1,*})]}$, and $\tilde{\sigma} = \sqrt{\mathbb{E}[h(\sigma^{1, 2, 1, *})h(\sigma^{1, 2, 2, *})]}$. Consider now the one-dimensional large portfolio model where the distance to default of the $i$-th asset $X_{t}^{i, *}$ evolves in time according to the system
\[
\begin{array}{rcl}
X_{t}^{i, *} &=& x^{i} + \left(r_i - \frac{\sigma_{2,1}^2}{2}\right)t + \tilde{\rho}_{1,i}\sigma_{2,1}W_{t}^{0} + \sqrt{1 - {\tilde{\rho}_{1,i}}^2}\sigma_{2,1}W_{t}^{i}, \,\, 0 \leq t \leq T_{i}^* \\
X_{t}^{i, *} &=& 0, \, t \geq T_{i}^*\\
T_{i}^* &=& \inf\{ t\geq 0: X_t^{i, *}=0\}, 
\end{array}
\]
where $\tilde{\rho}_{1,i} = \rho_{1,i}\frac{\tilde{\sigma}}{\sigma_{2,1}}$. Then, we have the convergence
\[
\mathbb{P}\Big[X_t^{1, \epsilon} \in \mathcal{I} \, \Big| \, W^{0}, \, B^{0}, \, \mathcal{G}\Big] \longrightarrow \mathbb{P}\Big[X_t^{1, *} \in \mathcal{I} \, \Big| \, W^{0}, \, \mathcal{G}\Big] \nonumber
\]
in distribution as $\epsilon \rightarrow 0^{+}$, for any interval $\mathcal{I} = \left(0, U\right]$ with $U \in \left(0, +\infty\right]$.
\end{thm}

\begin{rem}
Since all volatility processes have the same stationary distribution, a simple Cauchy-Schwartz inequality shows that $\tilde{\sigma} \leq \sigma_{2,1}$, which implies that $\tilde{\rho_{1,i}} \leq \rho_{1,i} < 1$ and $\sqrt{1 - {\tilde{\rho}_{1,i}}^2}$ is well-defined for each $i$.
\end{rem}

The above theorem gives only weak convergence and only under the restrictive assumption of having the same coefficients in each volatility. For this reason, we will also study the asymptotic behaviour of our system from a different perspective. In particular, we will fix the volatility path $\sigma^{1,1}$ and the coefficient vectors $C_i'$, and we will study the convergence of the solution $u^{\epsilon}(t, x)$ to the SPDE (\ref{1dimeq}) in 
the sped up setting, i.e
\begin{eqnarray}\label{eq:4.2}
u^{\epsilon}(t,\,x) &=& u_{0}(x) - \int_{0}^{t}\left(r-\frac{h^2\left(\sigma_{\frac{s}{\epsilon}}^{1,1}\right)}{2}\right)u_{x}^{\epsilon}(s,\,x)ds \nonumber \\
&& +\int_{0}^{t}\frac{h^2\left(\sigma_{\frac{s}{\epsilon}}^{1,1}\right)}{2}
u_{xx}^{\epsilon}(s,\,x)ds-\rho_{1,1}\int_{0}^{t}h\left(\sigma_{\frac{s}{\epsilon}}^{1,1}\right)u_{x}^{\epsilon}(s,\,x)dW_{s}^{0}, \nonumber \\
\end{eqnarray}
which is used to compute the loss $L_t^{\epsilon}$. 

We write now $\mathbb{E}_{\sigma, \mathcal{C}}$ to denote the expectation given the volatility path $\sigma^{1,1}$ and the $C_i$s, which we have fixed, and $L_{\sigma, \mathcal{C}}^2$ to denote the corresponding $L^2$ norms. By 2. of Theorem 4.1 in \cite{HK17}, the solution $u^{\epsilon}$ to the above SPDE satisfies the identity
\begin{equation}\label{eq:4.3}
\left\Vert u^{\epsilon}(t,\,\cdot)\right\Vert_{L^{2}(\mathbb{R}_{+})}^{2}+\left(1-\rho_{1,1}^{2}\right)\int_{0}^{t}h^2\left(\sigma_{\frac{t}{\epsilon}}^{1,1}\right)
\left\Vert u_{x}^{\epsilon}(s,\,\cdot)\right\Vert _{L^{2}(\mathbb{R}_{+})}^{2}ds=\left\Vert u_{0}\right\Vert _{L^{2}(\mathbb{R}_{+})}^{2},
\end{equation}
which shows that the $L^2(\mathbb{R}_{+})$ norms of $u^{\epsilon}$, and also its $L^2([0, \, T] \times \mathbb{R}_{+})$ norms (for any $T > 0$), are all uniformly bounded by a random variable which has a finite $L_{\sigma, \mathcal{C}}^2(\Omega)$ norm (the assumptions made in \cite{HK17} are also needed for this). It follows that in a subsequence of any given sequence of values of $\epsilon$ tending to zero, we have weak convergence to some element $u^{*}$ (see \cite{Brezis}), and we can have this both in $L_{\sigma, \mathcal{C}}^2([0, \, T] \times \mathbb{R}_{+} \times \Omega)$ and $\mathbb{P}$-almost surely in $L^2([0, \, T] \times \mathbb{R}_{+})$. The characterization of the weak limits $u^{*}$ is given in the following theorem.

\begin{thm}\label{thm1}
Suppose that $g$ has the positive recurrence property and that for some $C > 0$ we have $|h(x)| \leq C$ for all $x \geq 0$. Any weak limit $u^{*}$ of $u^{\epsilon}$ in $L_{\sigma, \mathcal{C}}^2([0, \, T] \times \mathbb{R}_{+} \times \Omega)$ solves the following SPDE 
\begin{eqnarray}\label{eq:lspde}
u^*(t,\,x) &=& u_{0}(x) - \left(r-\frac{\sigma_{2,1}^{2}}{2}\right)\int_{0}^{t}u_{x}^*(s,\,x)ds \nonumber \\
&& + \frac{\sigma_{2,1}^{2}}{2}\int_{0}^{t}
u_{xx}^*(s,\,x)ds-\rho_{1,1}\sigma_{1,1}\int_{0}^{t}u_{x}^*(s,\,x)dW_{s}^{0}.
\end{eqnarray}
If $h$ is bounded from below by a positive constant $c > 0$, the same weak convergence holds also in $H_{0}^{1}(\mathbb{R}_{+})\times L_{\sigma, \mathcal{C}}^2(\Omega \times [0, \, T])$, and $u^*$ is then the unique solution to (\ref{eq:lspde}) in that space. In this case there is a unique subsequential weak limit, and thus we have weak convergence as $\epsilon \rightarrow 0^{+}$. 
\end{thm}

It is not hard to see that the limiting SPDE (\ref{eq:lspde}) obtained in Theorem~\ref{thm1} corresponds to a constant volatility large portfolio model like the one given in Theorem~\ref{theorem} under the assumption that $(\kappa_{i},\,\theta_{i},\,v_{i},\,\rho_{2,i}) = (\kappa,\,\theta,\,v,\,\rho_{2})$, but with the correlation coefficients $\tilde{\rho}_{1,i} = \rho_{1,i}\frac{\tilde{\sigma}}{\sigma_{2,1}}$ replaced by $\rho_{1,i}' = \rho_{1,i}\frac{\sigma_{1,1}}{\sigma_{2,1}}$. This indicates that the convergence of the loss $L_t^{\epsilon}$ can only be established in a weak sense, as in general we will have $\tilde{\sigma} > \sigma_{1,1}$ and thus $\tilde{\rho}_{1,i} > \rho_{1,i}'$ for all $i$. This is stated explicitly in the next Proposition and its Corollary.

\begin{prop}\label{prp3}
Under the assumptions of Theorem~\ref{theorem}, we have always $\tilde{\sigma} \in [\sigma_{1,1}, \, \sigma_{2,1}]$. The lower and upper bounds are generally attained only when the volatilities are uncorrelated ($\rho_{2} = 0$) and perfectly correlated ($\rho_{2} \to 1$) respectively.
\end{prop}

\begin{cor}\label{cor1}
In general, the convergence established in Theorem~\ref{theorem} does not hold in any stronger sense, unless there is no market noise affecting all the volatilities in our setting.
\end{cor}

\section{The main results: small vol-of-vol setting}

We proceed now to the small vol-of-vol setting, where now only the volatility drifts are scaled by $\epsilon$, i.e $k_i = \kappa_i/\epsilon$ for all $i$. This leads to the model where the $i$-th asset's distance to default satisfies
\[
\begin{array}{rcl}
X_{t}^{i, \epsilon} &=& x^i + \int_{0}^{t}\left(r_{i}-\frac{h^{2}(\sigma_{t}^{i, \epsilon})}{2}\right)dt \nonumber \\
&& \qquad + \int_{0}^{t}h(\sigma_{t}^{i, \epsilon})\left(\sqrt{1-\rho_{1,i}^{2}}dW_{t}^{i}+\rho_{1,i}dW_{t}^{0}\right),\;\;0\leq t\leq T_{i}^{\epsilon} \\
\sigma_{t}^{i, \epsilon} &=& \sigma^{i,init} + \int_{0}^{t}\frac{\kappa_{i}}{\epsilon}(\theta_{i}-\sigma_{t}^{i, \epsilon})dt+\xi_{i}g\left(\sigma_{t}^{i, \epsilon}\right)\left(\sqrt{1-\rho_{2,i}^{2}}dB_{t}^{i}+\rho_{2,i}dB_{t}^{0}\right) \\
X_{t}^{i, \epsilon} &=& 0,\;t > T_{i}^{\epsilon} := \inf\{s \geq 0: X_s^{i, \epsilon} \leq 0\}. \\
\end{array}
\] 

The main feature of the above model is that when the random coefficients and the function $g$ satisfy certain conditions, the $i$-th volatility process $\sigma^{i, \epsilon}$ converges in a strong sense to the $\mathcal{C}$-measurable mean $\theta_i$ as $\epsilon \rightarrow 0^{+}$ for all $i \in \mathbb{N}$, and we can also determine the rate of convergence. The required conditions are the following, and they will be assumed to hold throughout the rest of this section:

\begin{enumerate}
    \item The i.i.d random variables $\sigma^{i}, \xi_i, \theta_i, \kappa_i$ take values in some compact subinterval of $\mathbb{R}$, with each $\kappa_i$ being bounded from below by some deterministic constant $c_{\kappa} > 0$.
    \item $g$ is a $C^1$ function with at most linear growth (i.e $\left|g(x)\right| \leq C_{1,g} + C_{2,g}|x|$ for some $C_{1,g}, C_{2,g} > 0$ and all $x \in \mathbb{R}$).
    \item Both the function $h$ and its derivative have polynomial growth.
\end{enumerate}

Under the above conditions, the convergence of each volatility process to its mean is given in the following proposition

\begin{prop}\label{lem1}
For any $t \geq 0$ and $p \geq 1$, we have $\sigma^{i, \epsilon} \rightarrow \theta_i$ as $\epsilon \rightarrow 0^{+}$ in $L^p(\Omega \times [0, \, t])$ at a rate of $\epsilon^{\frac{1}{p}}$. That is, we have $\Vert \sigma^{i, \epsilon} - \theta_i \Vert_{L^p(\Omega \times [0, \, t])}^p = \mathcal{O}(\epsilon)$ as $\epsilon \rightarrow 0^{+}$.
\end{prop} 
The reason for having only weak convergence of our system in the large vol-of-vol setting was the fact that the limiting quantities $\sigma_{1,1}$, $\sigma_{2,1}$ and $\tilde{\sigma}$ did not coincide. On the other hand, Proposition~\ref{lem1} implies that the corresponding limits in the small vol-of-vol setting are equal, allowing us to hope for our system to converge in a stronger sense. 

Let $u^{\epsilon}$ be the solution to the SPDE (\ref{1dimeq}) in the small vol-of-vol setting,
\begin{eqnarray}\label{eq:5.2}
u^{\epsilon}(t,\,x) &=& u_{0}(x) - \int_{0}^{t}\left(r-\frac{h^2\left(\sigma_s^{1,\epsilon}\right)}{2}\right)u_{x}^{\epsilon}(s,\,x)ds \nonumber \\
&& +\int_{0}^{t}\frac{h^2\left(\sigma_{s}^{1,\epsilon}\right)}{2}
u_{xx}^{\epsilon}(s,\,x)ds-\rho_{1,1}\int_{0}^{t}h\left(\sigma_{s}^{1,\epsilon}\right)u_{x}^{\epsilon}(s,\,x)dW_{s}^{0} \nonumber \\
\end{eqnarray}
where we have fixed the volatility paths and the random coefficients. Working as in the case of (\ref{eq:4.2}) and the proof of Theorem~\ref{thm4}, it is possible to establish similar asymptotic properties for the SPDE as $\epsilon \to 0^+$. However, it is more convenient to work with the antiderivative $v^{0, \epsilon} := \int_{\cdot}^{+\infty}u^{\epsilon}(\cdot, y)dy$, which satisfies the same SPDE but with different initial and boundary conditions, as the loss $L_t^{\epsilon} = 1 - \mathbb{P}[X_t^{1, \epsilon} > 0 \, | \, W^{0}, \, B^{0}, \, \mathcal{G}]$ equals the average of its value at $0$ over all possible volatility paths and coefficient values, while its convergence can be established in a much stronger sense and without the need to assume that $W^0$ and $B^0$ are uncorrelated. 
Our main result is stated below

\begin{thm}\label{thm5.2}
Define $v^{0}(t, x) = \int_{x}^{+\infty}u^0(t, y)dy$ for all $t, x \geq 0$ where $u^0$ is the unique solution to the SPDE 
\begin{eqnarray}\label{eq:5.2lim}
u^{0}(t,\,x) &=& u_{0}(x) - \int_{0}^{t}\left(r-\frac{h^2\left(\theta_1\right)}{2}\right)u_{x}^{0}(s,\,x)ds \nonumber \\
&& +\int_{0}^{t}\frac{h^2\left(\theta_1\right)}{2}
u_{xx}^{0}(s,\,x)ds-\rho_{1,1}\int_{0}^{t}h\left(\theta_1\right)u_{x}^{0}(s,\,x)dW_{s}^{0} \nonumber \\
\end{eqnarray}
in $L^2(\Omega \times [0, \, T] ; H_0^{1}(\mathbb{R}_{+}))$, which arises from the constant volatility model
\begin{equation}
\begin{array}{rcl}
dX_{t}^{i, *} &=& \left(r_{i}-\frac{h^{2}(\theta_i)}{2}\right)dt+h(\theta_i)\left(\sqrt{1-\rho_{1,i}^{2}}dW_{t}^{i}+\rho_{1,i}dW_{t}^{0}\right),\;\;0\leq t\leq T_{i} \\
X_{t}^{i, *} &=& 0,\;t > T_{i} := \inf\{s \geq 0: X_s^{i, *} \leq 0\} \\
X_{0}^{i, *} &=& x^{i}
\end{array}
\label{smallvisclim}
\end{equation} 
for $i \in \mathbb{N}$. Then, $v^{0, \epsilon}$ converges to $v^{0}$ as $\epsilon \rightarrow 0^{+}$, strongly in the Sobolev space $L^2(\Omega \times [0, \, T] ; H^{1}(\mathbb{R}_{+}))$ for any $T > 0$, and the rate of convergence is $\sqrt{\epsilon}$. That is, we have $\Vert v^{0, \epsilon} - v^{0} \Vert_{L^2(\Omega \times [0, \, T] ; H^{1}(\mathbb{R}_{+}))} = \mathcal{O}(\sqrt{\epsilon})$ as $\epsilon \to 0^+$
\end{thm}

The SPDE (\ref{eq:5.2lim}) corresponds to the model (\ref{smallvisclim}) in the sense that given the loss $L_t$, the mass of non-defaulted assets $1 - L_t$ equals
\[\mathbb{P}\Big[X_t^{1, *} > 0 \, \Big| \, W^{0}, \, \mathcal{G}\Big] = \mathbb{E}\Big[v^{0}\big(t, \, 0\big) \, \Big| \, W^0, \, \mathcal{G}\Big] = \int_{0}^{+\infty}\mathbb{E}\Big[u^{0}\big(t, \, x\big) \, \Big| \, W^0, \, \mathcal{G}\Big]dx.\]
In order to estimate the rate of convergence of probabilities of the form (\ref{lossprob}), we consider the approximation error
\begin{eqnarray}
&& E\big(x, \, T\big) = \nonumber \\
&& \int_{0}^{T}\Bigg|\mathbb{P}\bigg[\mathbb{P}\Big[X_t^{1, \epsilon} > 0 \, \Big| \, W^{0}, \, B^{0}, \, \mathcal{G}\Big] > x\bigg] - \mathbb{P}\bigg[\mathbb{P}\Big[X_t^{1, *} > 0 \, \Big| \, W^{0}, \, \mathcal{G}\Big] > x\bigg]\Bigg|dt \nonumber
\end{eqnarray}
for $x \in [0, \, 1]$, and determine its order of convergence. 

\begin{cor}\label{rateofcon}
For any $x \in [0, \, 1]$ such that 
$\mathbb{P}[X_t^{1, *} > 0 \, | \, W^{0}, \, C_{1}', \, \mathcal{G}]$ 
has a bounded density near $x$, uniformly in $t \in [0, \, T]$, we have 
$E(x, T) =\mathcal{O}(\epsilon^{1/3})$ as $\epsilon \to 0^{+}$
\end{cor}

\section{Proofs: large vol-of-vol setting}
We prove Theorem~\ref{theorem}, Theorem~\ref{thm1}, Proposition~\ref{prp3} and Corollary~\ref{cor1}, the main results of section 2.

\begin{proof}[\textbf{Theorem 2.4}]
To establish convergence in distribution, we show that, for every bounded and continuous function $G: \mathbb{R} \rightarrow \mathbb{R}$, we have:
\begin{eqnarray}\label{distconvmain}
&& \mathbb{E}\Bigg[G\bigg(\mathbb{P}\Big[X_t^{1, \epsilon} \in \mathcal{I} \, \Big| \, W^{0}, \, B^{0}, \, \mathcal{G}\Big]\bigg) \Bigg] \longrightarrow \mathbb{E}\Bigg[G\bigg(\mathbb{P}\Big[X_t^{1, *} \in \mathcal{I} \, \Big| \, W^{0}, \, \mathcal{G}\Big]\bigg) \Bigg] \nonumber \\
\end{eqnarray} as $\epsilon \rightarrow 0^{+}$, where $\mathcal{I} = (0, U]$. Observe now that since the conditional probabilities take values in the compact interval $[0, 1]$, it is equivalent to have (\ref{distconvmain}) for all continuous $G: [0, 1] \rightarrow \mathbb{R}$, and by the Weierstrass approximation theorem and linearity, we actually need to have this only when $G$ is a polynomial of the form $G(x)=x^m$.
We now write $Y^{i, \epsilon}$ for the $i$-th asset's distance to default in the sped up volatility setting, when the stopping condition at zero is ignored, that is
\begin{eqnarray}
Y_{t}^{i, \epsilon} &=& x^{i} + \int_{0}^{t}\left(r_{i}-\frac{h^{2}\left(\sigma_{\frac{s}{\epsilon}}^{i,1}\right)}{2}\right)ds \nonumber \\
&& +\int_{0}^{t}h(\sigma_{\frac{s}{\epsilon}}^{i,1})\rho_{1,i}dW_{s}^{0} +\int_{0}^{t}h(\sigma_{\frac{s}{\epsilon}}^{i,1})\sqrt{1-\rho_{1,i}^{2}}dW_{s}^{i} \nonumber
\end{eqnarray}
with
\begin{eqnarray}
\sigma_{t}^{i,1} &=& \sigma^{i,init} + \kappa\int_{0}^{t}\left(\theta - \sigma_{s}^{i,1}\right)ds \nonumber \\
&& \qquad + v\int_{0}^{t}g\left(\sigma_{s}^{i,1}\right)\rho_2dB_{s}^{0} + v\int_{0}^{t}g\left(\sigma_{s}^{i,1}\right)\sqrt{1 - \rho_2^2}dB_{s}^{i} \nonumber
\end{eqnarray}
for all $t \geq 0$, and then we have $X_{t}^{i, \epsilon} =  Y_{t \wedge T_i^{\epsilon}}^{i, \epsilon}$. The $m$ stochastic processes $\{X^{i, \epsilon}: \, 1 \leq i \leq m \}$ are obviously pairwise i.i.d when the information contained in $W^{0}, B^{0}$ and $\mathcal{G}$ is given. Therefore we can write:
\begin{eqnarray}\label{rep1}
&& \mathbb{E}\Bigg[G\bigg(\mathbb{P}\Big[X_t^{1, \epsilon} \in \mathcal{I} \, \Big| \, W^{0}, \, B^{0}, \, \mathcal{G}\Big]\bigg)\Bigg]  \nonumber \\
&& \qquad \qquad = \mathbb{E}\bigg[\mathbb{P}^m\Big[X_t^{1, \epsilon} \in \mathcal{I} \, \Big| \, W^{0}, \, B^{0}, \, \mathcal{G}\Big] \bigg]  \nonumber \\
&& \qquad \qquad = \mathbb{E}\bigg[\mathbb{P}\Big[X_t^{1, \epsilon} \in \mathcal{I}, \, X_t^{2, \epsilon} \in \mathcal{I}, \, ..., \, X_t^{m, \epsilon} \in \mathcal{I}  \, \Big| \, W^{0}, \, B^{0}, \, \mathcal{G}\Big] \bigg]  \nonumber \\
&& \qquad \qquad = \mathbb{P}\Big[X_t^{1, \epsilon} \in \mathcal{I}, \, X_t^{2, \epsilon} \in \mathcal{I}, \, ..., \, X_t^{m, \epsilon} \in \mathcal{I}  \Big]  \nonumber \\
&& \qquad \qquad = \mathbb{P}\bigg[\Big(\min_{1 \leq i \leq m}\min_{0 \leq s \leq t}Y_s^{i, \epsilon} \, , \max_{1 \leq i \leq m}Y_t^{i, \epsilon} \Big) \in \left(0, \, +\infty\right)\times \left(-\infty, U\right] \bigg]. \nonumber \\
\end{eqnarray}
Next, for each $i$, we write $Y^{i, *}$ for the process $X^{i, *}$ when the stopping condition at zero is ignored, that is
\begin{eqnarray}
Y_{t}^{i, *} = X_0^{i} + \left(r_i - \frac{\sigma_{2,1}^2}{2}\right)t + \tilde{\rho}_{1,i}\sigma_{2,1}W_{t}^{0} + \sqrt{1 - {\tilde{\rho}_{1,i}}^2}\sigma_{2,1}W_{t}^{i} \nonumber
\end{eqnarray}
for all $t \geq 0$, with $\tilde{\rho}_{1,i} = \rho_{1,i}\frac{\tilde{\sigma}}{\sigma_{2,1}}$. Again, it is easy to check that the processes $Y^{i,*}$ are pairwise i.i.d when the information contained in $W^{0}, B^{0}$ and $\mathcal{G}$ is given. Thus, we can write
\begin{eqnarray}\label{rep2}
&& \mathbb{E}\Bigg[G\bigg(\mathbb{P}\Big[X_t^{1, *} \in \mathcal{I} \, \Big| \, W^{0}, \, \mathcal{G}\Big]\bigg)\Bigg]  \nonumber \\
&& \qquad \qquad = \mathbb{E}\bigg[\mathbb{P}^m\Big[X_t^{1, *} \in \mathcal{I} \, \Big| \, W^{0}, \, \mathcal{G}\Big] \bigg]  \nonumber \\
&& \qquad \qquad = \mathbb{E}\bigg[\mathbb{P}\Big[X_t^{1, *} \in \mathcal{I}, \, X_t^{2, *} \in \mathcal{I}, \, ..., \, X_t^{m, *} \in \mathcal{I}  \, \Big| \, W^{0}, \, \mathcal{G}\Big] \bigg]  \nonumber \\
&& \qquad \qquad = \mathbb{P}\Big[X_t^{1, *} \in \mathcal{I}, \, X_t^{2, *} \in \mathcal{I}, \, ..., \, X_t^{m, *} \in \mathcal{I}  \Big]  \nonumber \\
&& \qquad \qquad = \mathbb{P}\bigg[\Big(\min_{1 \leq i \leq m}\min_{0 \leq s \leq t}Y_s^{i, *} \, , \max_{1 \leq i \leq m}Y_t^{i, *} \Big) \in \left(0, \, +\infty\right)\times \left(-\infty, U\right] \bigg]. \nonumber \\
\end{eqnarray}
Then, (\ref{rep1}) and (\ref{rep2}) show that the result we want to prove has been reduced to the convergence 
\[\left(\min_{1 \leq i \leq m}\min_{0 \leq s \leq t}Y_s^{i, \epsilon} \, , \max_{1 \leq i \leq m}Y_t^{i, \epsilon} \right) \longrightarrow \left(\min_{1 \leq i \leq m}\min_{0 \leq s \leq t}Y_s^{i, *} \, , \max_{1 \leq i \leq m}Y_t^{i, *} \right),\]
in distribution as $\epsilon \rightarrow 0^{+}$ (since the probability that any of the $m$ minimums equals zero is zero, as the minimum of any Gaussian process is always continuously distributed, while $Y^{i, \epsilon}$ is obviously Gaussian for any given path of $\sigma^{i,1}$).

Let $C([0, \, t]; \mathbb{R}^{m})$ be the classical Wiener space of continuous functions defined on $[0, \, t]$ and taking values in $\mathbb{R}^{m}$ (i.e the space of these functions equipped with the supremum norm and the Wiener probability measure), and observe that $\min_{1 \leq i \leq m}p_i(\min_{0 \leq s \leq t}\cdot (s))$ defined on $C([0, \, t]; \mathbb{R}^{m})$, where $p_i$ stands for the projection on the $i$-th axis, is a continuous functional. Indeed, for any two continuous functions $f_1, \, f_2$ in $C([0, \, t]; \mathbb{R}^{m})$, we have: 
\begin{eqnarray}
&& \Bigg|\min_{1 \leq i \leq m}p_i\bigg(\min_{0 \leq s \leq t}f_1(s)\bigg) - \min_{1 \leq i \leq m}p_i\bigg(\min_{0 \leq s \leq t}f_2(s)\bigg)\Bigg| \nonumber \\
&& \qquad \quad \qquad \qquad \qquad = \Big|p_{i_1}\big(f_1(s_1)\big) - p_{i_2}\big(f_2(s_2)\big)\Big| \nonumber
\end{eqnarray}
for some $s_1, s_2 \in [0, \, t]$ and $1 \leq i_1, i_2 \leq m$, and without loss of generality we may assume that the difference inside the last absolute value is nonnegative. Moreover we have:
\begin{eqnarray}
p_{i_1}\big(f_1(s_1)\big) = \min_{1 \leq i \leq m}p_i\bigg(\min_{0 \leq s \leq t}f_1(s)\bigg) \leq p_{i_2}\big(f_1(s_2)\big) \nonumber
\end{eqnarray}
and thus 
\begin{eqnarray}
&& \Bigg|\min_{1 \leq i \leq m}p_i\bigg(\min_{0 \leq s \leq t}f_1(s)\bigg) - \min_{1 \leq i \leq m}p_i\bigg(\min_{0 \leq s \leq t}f_2(s)\bigg)\Bigg| \nonumber \\
&& \quad \qquad \qquad \qquad \qquad = p_{i_1}\big(f_1(s_1)\big) - p_{i_2}\big(f_2(s_2)\big)  \nonumber \\
&& \quad \qquad \qquad \qquad \qquad \leq p_{i_2}\big(f_1(s_2)\big) - p_{i_2}\big(f_2(s_2)\big) \nonumber \\
&& \quad \qquad \qquad \qquad \qquad \leq \Big|p_{i_2}\big(f_1(s_2)\big) - p_{i_2}\big(f_2(s_2)\big)\Big| \nonumber \\
&& \quad \qquad \qquad \qquad \qquad \leq \left\Vert f_1 - f_2 \right\Vert_{C\left(\left[0, \, t\right]; \mathbb{R}^{m}\right)}. \nonumber
\end{eqnarray}
Obviously, $\max_{1 \leq i \leq m}p_i(\cdot (t))$ defined on $C([0, \, t]; \mathbb{R}^{m})$ is also continuous (as the maximum of finitely many evaluation functionals).
Therefore, our problem is finally reduced to showing that $(Y^{1, \epsilon}, \, Y^{2, \epsilon}, \, ..., \, Y^{m, \epsilon})$ converges in distribution to $(Y^{1, *}, \, Y^{2, *}, \, ..., \, Y^{m, *})$ in the space $C([0, \, t]; \mathbb{R}^{m})$, as $\epsilon \rightarrow 0^{+}$. 

In order to show the convergence in distribution we first establish that a limit in distribution exists as as $\epsilon \rightarrow 0^{+}$ by using a tightness argument, and then we will characterize the limits of the finite dimensional distributions. To show tightness of the laws of $(Y^{1, \epsilon}, \, Y^{2, \epsilon}, \, ..., \, Y^{m, \epsilon})$ for $\epsilon \in \mathbb{R}_{+}$, which implies the desired convergence in distribution, we recall a special case of Theorem 3.7.2 in Ethier and Kurtz \cite{EK05} for continuous processes, according to which it suffices to prove that for a given $\eta > 0$, there exist some $\delta > 0$ and $N > 0$ such that:
\begin{eqnarray}\label{cond1}
\mathbb{P}\Bigg[\bigg| \Big(Y_{0}^{1, \epsilon}, \, Y_{0}^{2, \epsilon}, \, ..., \, Y_{0}^{m, \epsilon}\Big) \bigg|_{\mathbb{R}^{m}} > N \Bigg] \leq \eta 
\end{eqnarray}
and
\begin{eqnarray}\label{cond2}
&& \mathbb{P}\Bigg[\sup_{0\leq s_1, s_2 \leq t, \, |s_1 - s_2| \leq \delta}\bigg| \Big(Y_{s_1}^{1, \epsilon}, \, Y_{s_1}^{2, \epsilon}, \, ..., \, Y_{s_1}^{m, \epsilon}\Big) \nonumber \\
&& \qquad \qquad \qquad \qquad \qquad \quad - \Big(Y_{s_2}^{1, \epsilon}, \, Y_{s_2}^{2, \epsilon}, \, ..., \, Y_{s_2}^{m, \epsilon} \Big)\bigg|_{\mathbb{R}^{m}} > \eta \Bigg] \leq \eta
\end{eqnarray}
for all $\epsilon > 0$. (\ref{cond1}) can easily be achieved for some very large $N > 0$, since we have $(Y_{0}^{1, \epsilon}, \, Y_{0}^{2, \epsilon}, \, ..., \, Y_{0}^{m, \epsilon}) = (x^{1}, \, x^{2}, \, ..., \, x^{m})$, which is independent of $\epsilon$ and almost surely finite (the sum of the probabilities that the norm of this vector belongs to $\left[n, \, n+1 \right]$ over $n \in \mathbb{N}$ is a convergent series and thus, by the Cauchy criteria, the same sum but for $n \geq N$ tends to zero as $N$ tends to infinity). For (\ref{cond2}), observe that $| \cdot |_{\mathbb{R}^{m}}$ can be any of the standard equivalent $L^p$ norms of $\mathbb{R}^{m}$, and we choose it to be $L^{\infty}$. Then we have: 
\begin{eqnarray}\label{firstestimate4}
&& \mathbb{P}\Bigg[\sup_{0\leq s_1, s_2 \leq t, \, |s_1 - s_2| \leq \delta}\Big| \big(Y_{s_1}^{1, \epsilon}, \, Y_{s_1}^{2, \epsilon}, \, ..., \, Y_{s_1}^{m, \epsilon}\big) \nonumber \\
&& \qquad \qquad \qquad \qquad \quad \qquad - \big(Y_{s_2}^{1, \epsilon}, \, Y_{s_2}^{2, \epsilon}, \, ..., \, Y_{s_2}^{m, \epsilon} \big)\Big|_{\mathbb{R}^{m}} > \eta \Bigg] \nonumber \\
&& \qquad \qquad \qquad = \mathbb{P}\Bigg[\cup_{i=1}^{m}\bigg\{\sup_{0\leq s_1, s_2 \leq t, \, |s_1 - s_2| \leq \delta}\big|Y_{s_1}^{i, \epsilon} - Y_{s_2}^{i, \epsilon}\big| > \eta\bigg\}  \Bigg] \nonumber \\ 
&& \qquad \qquad \qquad \leq \sum_{i=1}^{m}\mathbb{P}\Bigg[\sup_{0\leq s_1, s_2 \leq t, \, |s_1 - s_2| \leq \delta}\big|Y_{s_1}^{i, \epsilon} - Y_{s_2}^{i, \epsilon}\big| > \eta  \Bigg] \nonumber \\ 
&& \qquad \qquad \qquad = m\mathbb{P}\Bigg[\sup_{0\leq s_1, s_2 \leq t, \, |s_1 - s_2| \leq \delta}\left|Y_{s_1}^{1, \epsilon} - Y_{s_2}^{1, \epsilon}\right| > \eta  \Bigg]  
\end{eqnarray}
and since the Ito integral $\int_{0}^{t}h(\sigma_{\frac{s}{\epsilon}}^{1, 1})d(\sqrt{1-\rho_{1}^{2}}W_{s}^{1}+\rho_{1}W_{s}^{0})$ can be written as $\tilde{W}_{\int_{0}^{t}h^2(\sigma_{\frac{s}{\epsilon}}^{1,1})ds}$, where $\tilde{W}$ is another standard Brownian motion, denoting the maximum of $h$ by $M$ we also have:
\begin{eqnarray*}
&& \mathbb{P}\left[\sup_{0\leq s_1, s_2 \leq t, \, |s_1 - s_2| \leq \delta}\left|Y_{s_1}^{1, \epsilon} - Y_{s_2}^{1, \epsilon}\right| > \eta  \right] \nonumber \\
&& \qquad = \mathbb{P}\Bigg[\sup_{0\leq s_1, s_2 \leq t, \, |s_1 - s_2| \leq \delta}\Bigg|\int_{s_2}^{s_1}\left(r-\frac{h^{2}\left(\sigma_{\frac{s}{\epsilon}}^{1,1}\right)}{2}\right)ds \nonumber \\
&& \qquad \qquad \qquad \qquad \qquad \qquad \quad +\left(\tilde{W}_{\int_{0}^{s_1}h^2\left(\sigma_{\frac{s}{\epsilon}}^{1,1}\right)ds} - \tilde{W}_{\int_{0}^{s_2}h^2\left(\sigma_{\frac{s}{\epsilon}}^{1,1}\right)}\right)\Bigg| > \eta  \Bigg] \nonumber \\
&& \qquad \leq \mathbb{P}\left[\sup_{0\leq s_1, s_2 \leq t, \, |s_1 - s_2| \leq \delta}\left|\int_{s_2}^{s_1}\left(r-\frac{h^{2}\left(\sigma_{\frac{s}{\epsilon}}^{1,1}\right)}{2}\right)ds\right| > \frac{\eta}{2} \right] \nonumber \\
&& \qquad \qquad + \mathbb{P}\Bigg[\sup_{0\leq s_1, s_2 \leq t, \, |s_1 - s_2| \leq \delta}\Bigg|\tilde{W}_{\int_{0}^{s_1}h^2\left(\sigma_{\frac{s}{\epsilon}}^{1,1}\right)ds} - \tilde{W}_{\int_{0}^{s_2}h^2\left(\sigma_{\frac{s}{\epsilon}}^{1,1}\right)ds}\Bigg| > \frac{\eta}{2} \Bigg] \nonumber \\
&& \qquad \leq \mathbb{P}\left[\delta\left(r + M\right) > \frac{\eta}{2} \right] \nonumber \\
&& \qquad \qquad + \mathbb{P}\left[\sup_{0\leq s_3, s_4 \leq M^2t, \, |s_3 - s_4| \leq M^2\delta}\left|\tilde{W}_{s_{3}} - \tilde{W}_{s_{4}}\right| > \frac{\eta}{2} \right] 
\end{eqnarray*}
since $|\int_{a}^{b}h^2(\sigma_{\frac{s}{\epsilon}}^{1,1})ds| \leq M^2|a - b|$ for all $a, b \in \mathbb{R}_{+}$. The first of the last two probabilities is clearly zero for $\delta < \frac{\eta}{2(r + M)}$, while the second one can also be made arbitrarily small for small enough $\delta$, since by a well known result about the modulus of continuity of a Brownian motion (see Levy \cite{PLEV}) the supremum within that probability converges almost surely (and thus also in probability) to $0$ as fast as $M\sqrt{2\delta \ln \frac{1}{M^2\delta}}$. Using these in (\ref{firstestimate4}) we deduce that (\ref{cond2}) is also satisfied and we have the desired tightness result, which implies that $(Y_{\cdot}^{1, \epsilon}, \, ..., \, Y_{\cdot}^{m, \epsilon})$ converges in distribution to some limit $(Y_{\cdot}^{1, 0}, \, ..., \, Y_{\cdot}^{m, 0})$ (along some sequence). 

To conclude our proof, we need to show that $(Y_{\cdot}^{1, 0}, \, ..., \, Y_{\cdot}^{m, 0})$ and $(Y_{\cdot}^{1, *}, \, ..., \, Y_{\cdot}^{m, *})$ coincide. Since both $m$-dimensional processes are uniquely determined by their finite-dimensional distributions, and since evaluation functionals on $C([0, \, t]; \mathbb{R}^{m})$ preserve convergences in distribution (as continuous linear functionals), we only need to show that for any fixed $(i_1, \, ..., \, i_{\ell}) \in \{1, \, ..., \, m\}^{\ell}$, any fixed $(t_1, \, ..., \, t_{\ell}) \in (0, \, +\infty)^{\ell}$, and any fixed continuous and bounded function $q: \mathbb{R}^{\ell} \rightarrow \mathbb{R}$, for an arbitrary $\ell \in \mathbb{N}$, we have
\begin{eqnarray}
&& \mathbb{E}\bigg[q\left(Y_{t_1}^{i_1, \epsilon}, \, Y_{t_2}^{i_2, \epsilon}, \, ..., \, Y_{t_{\ell}}^{i_{\ell}, \epsilon}\right)\bigg] \longrightarrow \mathbb{E}\bigg[q\left(Y_{t_1}^{i_1, *}, \, Y_{t_2}^{i_2, *}, \, ..., \, Y_{t_{\ell}}^{i_{\ell}, *}\right)\bigg] \nonumber
\end{eqnarray}
as $\epsilon \rightarrow 0^{+}$. By the dominated convergence theorem, the above follows if we are able to show that 
\begin{eqnarray}
&& \lim_{\epsilon \rightarrow 0^{+}}\mathbb{E}\bigg[q\left(Y_{t_1}^{i_1, \epsilon}, \, Y_{t_2}^{i_2, \epsilon}, \, ..., \, Y_{t_{\ell}}^{i_{\ell}, \epsilon}\right) \,\bigg| \sigma_{\cdot}^{i_1, 1}, \, \sigma_{\cdot}^{i_2, 1}, \, ..., \, \sigma_{\cdot}^{i_{\ell}, 1}, \, \mathcal{C} \bigg] \nonumber \\
&& \qquad = \mathbb{E}\bigg[q\left(Y_{t_1}^{i_1, *}, \, Y_{t_2}^{i_2, *}, \, ..., \, Y_{t_{\ell}}^{i_{\ell}, *}\right) \,\bigg| \sigma_{\cdot}^{i_1, 1}, \, \sigma_{\cdot}^{i_2, 1}, \, ..., \, \sigma_{\cdot}^{i_{\ell}, 1}, \, \mathcal{C} \bigg] \nonumber
\end{eqnarray}
$\mathbb{P}$- almost surely. However, when the information contained in $\sigma_{\cdot}^{i_1, 1}, \, ..., \, \sigma_{\cdot}^{i_{\ell}, 1}$ and $\mathcal{C}$ is given, both $(Y_{t_1}^{i_1, \epsilon}, \, ..., \, Y_{t_{\ell}}^{i_{\ell}, \epsilon})$ and $(Y_{t_1}^{i_1, *}, \, ..., \, Y_{t_{\ell}}^{i_{\ell}, *})$ follow a normal distribution in $\mathbb{R}^{\ell}$. This means that given $(\sigma_{\cdot}^{i_1, 1}, \, ..., \, \sigma_{\cdot}^{i_{\ell}, 1})$ and $\mathcal{C}$, we only need to show that as $\epsilon \rightarrow 0^{+}$, the mean vector and the covariance matrix of $(Y_{t_1}^{i_1, \epsilon}, \, ..., \, Y_{t_{\ell}}^{i_{\ell}, \epsilon})$ converge to the mean vector and the covariance matrix of $(Y_{t_1}^{i_1, *}, \, ..., \, Y_{t_{\ell}}^{i_{\ell}, *})$ respectively. Given $(\sigma_{\cdot}^{i_1, 1}, \, ..., \, \sigma_{\cdot}^{i_{\ell}, 1})$, the information contained in $\mathcal{C}$, and a $k \in \{1, \, 2, \, ..., \, \ell\}$, the $k$-th coordinate of the mean vector of $(Y_{t_1}^{i_1, \epsilon}, \, Y_{t_2}^{i_2, \epsilon}, \, ..., \, Y_{t_{\ell}}^{i_{\ell}, \epsilon})$ is equal to $X_{0}^{i_k} + \int_{0}^{t_k}(r_{i_k}-\frac{h^{2}(\sigma_{\frac{s}{\epsilon}}^{i_k, 1})}{2})ds$, and by the positive recurrence property it converges as $\epsilon \rightarrow 0^{+}$ to $X_0^{i_k} + (r_{i_k} - \frac{\sigma_{2,1}^2}{2})t_k$ (since the volatility processes all have the same coefficients and thus the same stationary distributions), which is the $k$-th coordinate of the mean vector of $(Y_{t_1}^{i_1, *}, \, Y_{t_2}^{i_2, *}, \, ..., \, Y_{t_{\ell}}^{i_{\ell}, *})$. Now we only need to obtain the corresponding convergence result for the covariance matrices of our processes. For some $1 \leq p, q \leq \ell$, given $(\sigma_{\cdot}^{i_1, 1}, \, \sigma_{\cdot}^{i_2, 1}, \, ..., \, \sigma_{\cdot}^{i_{\ell}, 1})$ and the information contained in $\mathcal{C}$, the covariance of $Y_{t_p}^{i_p, \epsilon}$ and $Y_{t_q}^{i_q, \epsilon}$ is equal to
\begin{eqnarray}
&& \left(\rho_{1, i_p}\rho_{1, i_q} + \delta_{i_p, i_q}\sqrt{1 - \rho_{1, i_p}}\sqrt{1 - \rho_{1, i_q}}\right)\int_{0}^{t_p \wedge t_q}h\left(\sigma_{\frac{s}{\epsilon}}^{i_p, 1}\right)h\left(\sigma_{\frac{s}{\epsilon}}^{i_q, 1}\right)ds, \nonumber
\end{eqnarray} 
while the covariance of $Y_{t_p}^{i_p, *}$ and $Y_{t_q}^{i_q, *}$ is equal to
\begin{eqnarray}
&& \left(\tilde{\rho}_{1,i_p}\tilde{\rho}_{1,i_q} + \delta_{i_p, i_q}\sqrt{1 - \tilde{\rho}_{1, i_p}^2}\sqrt{1 - \tilde{\rho}_{1, i_q}^2}\right) \sigma_{2, 1}^2 t_p \wedge t_q. \nonumber
\end{eqnarray}
This means that for $i_p = i_q = i \in \{1, \, 2, \, ..., \, m\}$ we need to show that
\begin{eqnarray}
&& \int_{0}^{t_p \wedge t_q}h^2\left(\sigma_{\frac{s}{\epsilon}}^{i, 1}\right)ds \longrightarrow \sigma_{2, 1}^2 t_p \wedge t_q \nonumber
\end{eqnarray}
as $\epsilon \rightarrow 0^{+}$, while for $i_p \neq i_q$ we need to show that:
\begin{eqnarray}
&& \rho_{1, i_p}\rho_{1, i_q} \int_{0}^{t_p \wedge t_q}h\left(\sigma_{\frac{s}{\epsilon}}^{i_p, 1}\right)h\left(\sigma_{\frac{s}{\epsilon}}^{i_q, 1}\right)ds \longrightarrow \tilde{\rho}_{1, i_p}\tilde{\rho}_{1, i_q} \sigma_{2, 1}^2 t_p \wedge t_q \nonumber
\end{eqnarray}
as $\epsilon \rightarrow 0^{+}$, where $\tilde{\rho}_{1, i}\sigma_{2, 1} = \rho_{1, i}\tilde{\sigma}$ for all $i \leq m$. Both convergence results follow from the positive recurrence property for $\tilde{\sigma} = \sqrt{\mathbb{E}\left[h\left(\sigma^{i_p, i_q, 1, *}\right)h\left(\sigma^{i_p, i_q, 2, *}\right)\right]}$, which does not depend on $i_p$ and $i_q$ since the volatility processes all have the same coefficients and thus the same joint stationary distributions. This concludes the proof.
\qed
\end{proof}

\begin{proof}[\textbf{Theorem 2.6}]
Let $\mathbb{V}$ be the set of $W_{\cdot}^0$-adapted, square-integrable semimartingales on $[0, \, T]$. Thus for any $\{V_t: \, 0 \leq t \leq T\} \in \mathbb{V}$, there exist two $W^0$-adapted and square-integrable processes $\{v_{1,t}: \, 0 \leq t \leq T\}$ and $\{v_{2,t}: \, 0 \leq t \leq T\}$, such that
\begin{eqnarray}\label{d1}
V_t = V_0 + \int_{0}^{t}v_{1,s}ds + \int_{0}^{t}v_{2,s}dW_{s}^{0},
\end{eqnarray}
for all $t \geq 0$. The processes of the above form for which $\{v_{1,t}: \, 0 \leq t \leq T\}$ and $\{v_{2,t}: \, 0 \leq t \leq T\}$ are simple processes, that is
\begin{eqnarray}\label{d2}
v_{i,t} = F_i\mathbb{I}_{\left[t_{1}, \, t_{2}\right]}(t),
\end{eqnarray}
for all $0 \leq t \leq T$ and $i \in \{1, \, 2\}$, with each $F_i$ being $\mathcal{F}_{t_{1}}^{W^0}$-measurable, span a linear subspace $\tilde{\mathbb{V}}$ which is dense in $\mathbb{V}$ under the $L^2$ norm. By using the boundedness of $h$ and then the estimate (\ref{eq:4.3}), for any $p > 0$ and any $T > 0$ we obtain
\begin{eqnarray}\label{bound}
&& \int_{0}^{T}\bigg\Vert h^{p}\Big(\sigma_{\frac{t}{\epsilon}}^{1,1}\Big)u^{\epsilon}(t,\,\cdot)\bigg\Vert_{L_{\sigma, \mathcal{C}}^{2}(\mathbb{R}_{+}\times \Omega)}^{2}dt \leq TC^{2p}\left\Vert u_{0}\right\Vert _{L^{2}(\mathbb{R}_{+})}^{2}.
\end{eqnarray}
It follows that given a sequence $\epsilon_{n} \rightarrow 0^{+}$, there exists always a subsequence $\{\epsilon_{k_{n}}: \, n \in \mathbb{N}\}$, such that $h^{p}(\sigma_{\frac{\cdot}{\epsilon}}^{1,1})u^{\epsilon}(\cdot,\,\cdot)$ converges weakly to some $u_{p}(\cdot, \cdot)$ in the space $L_{\sigma, \mathcal{C}}^2([0, \, T] \times \mathbb{R}_{+}\times \Omega)$ for $p \in \{1, \, 2 \}$. Testing (\ref{eq:4.2}) against an arbitrary smooth and compactly supported function $f$ of $x \in \mathbb{R}_{+}$, using Ito's formula for the product of $\int_{\mathbb{R}_{+}}u^{\epsilon}(\cdot, \, x)f(x)dx$ with a process $V_{\cdot} \in \tilde{\mathbb{V}}$ having the form (\ref{d1}) - (\ref{d2}), and finally taking expectations, we find that: 
\begin{eqnarray}\label{tested}
&& \mathbb{E}_{\sigma, \mathcal{C}}\left[V_t\int_{\mathbb{R}_{+}}u^{\epsilon}(t,\,x)f(x)dx\right] \nonumber \\
&& \qquad = \mathbb{E}_{\sigma, \mathcal{C}}\left[V_0\int_{\mathbb{R}_{+}}u_{0}(x)f(x)dx\right] + r\int_{0}^{t}\mathbb{E}_{\sigma, \mathcal{C}}\left[V_s\int_{\mathbb{R}_{+}}u^{\epsilon}(s,\,x)f'(x)dx\right]ds \nonumber \\
&& \qquad \qquad - \int_{0}^{t}\mathbb{E}_{\sigma, \mathcal{C}}\left[V_s\int_{\mathbb{R}_{+}}\frac{h^2\left(\sigma_{\frac{s}{\epsilon}}^{1,1}\right)}{2}u^{\epsilon}(s,\,x)f'(x)dx\right]ds \nonumber \\
&& \qquad \qquad +\int_{0}^{t}\mathbb{E}_{\sigma, \mathcal{C}}\left[V_s\int_{\mathbb{R}_{+}}\frac{h^2\left(\sigma_{\frac{s}{\epsilon}}^{1,1}\right)}{2}
u^{\epsilon}(s,\,x)f''(x)dx\right]ds \nonumber \\
&& \qquad \qquad +\int_{0}^{t}\mathbb{E}_{\sigma, \mathcal{C}}\left[v_{1,s}\int_{\mathbb{R}_{+}}
u^{\epsilon}(s,\,x)f(x)dx\right]ds \nonumber \\
&& \qquad \qquad + \rho_{1,1}\int_{0}^{t}\mathbb{E}_{\sigma, \mathcal{C}}\left[v_{2,s}\int_{\mathbb{R}_{+}}h\left(\sigma_{\frac{s}{\epsilon}}^{1,1}\right)u^{\epsilon}(s,\,x)f'(x)dx\right]ds
\end{eqnarray}
for all $t \leq T$. Thus, setting $\epsilon = \epsilon_{k_{n}}$ and taking $n \rightarrow +\infty$, by the weak convergence results mentioned above we obtain
\begin{eqnarray}\label{testedlimit}
&& \mathbb{E}_{\sigma, \mathcal{C}}\left[V_t\int_{\mathbb{R}_{+}}u^*(t,\,x)f(x)dx\right] \nonumber \\
&& \qquad = \mathbb{E}_{\sigma, \mathcal{C}}\left[V_0\int_{\mathbb{R}_{+}}u_{0}(x)f(x)dx\right] + r\int_{0}^{t}\mathbb{E}_{\sigma, \mathcal{C}}\left[V_s\int_{\mathbb{R}_{+}}u^*(s,\,x)f'(x)dx\right]ds \nonumber \\
&& \qquad \qquad - \frac{1}{2}\int_{0}^{t}\mathbb{E}_{\sigma, \mathcal{C}}\left[V_s\int_{\mathbb{R}_{+}}u_{2}(s,\,x)f'(x)dx\right]ds \nonumber \\
&& \qquad \qquad + \frac{1}{2}\int_{0}^{t}\mathbb{E}_{\sigma, \mathcal{C}}\left[V_s\int_{\mathbb{R}_{+}}
u_{2}(s,\,x)f''(x)dx\right]ds \nonumber \\
&& \qquad \qquad +\int_{0}^{t}\mathbb{E}_{\sigma, \mathcal{C}}\left[v_{1,s}\int_{\mathbb{R}_{+}}
u^*(s,\,x)f(x)dx\right]ds \nonumber \\
&& \qquad \qquad + \rho_{1,1}\int_{0}^{t}\mathbb{E}_{\sigma, \mathcal{C}}\left[v_{2,s}\int_{\mathbb{R}_{+}}u_{1}(s,\,x)f'(x)dx\right]ds
\end{eqnarray}
for all $0 \leq t \leq T$. The convergence of the terms in the RHS of (\ref{tested}) holds pointwise in $t$, while the one term in the LHS converges weakly. Since we can easily find uniform bounds for all the terms in (\ref{tested}) (by using (\ref{bound})), the dominated convergence theorem implies that all the weak limits coincide with the corresponding pointwise limits, which gives (\ref{testedlimit}) as a limit of (\ref{tested}) both weakly and pointwise in $t$. It is clear then that $\mathbb{E}_{\sigma, \mathcal{C}}[V_t\int_{\mathbb{R}_{+}}u^*(t,\,x)f(x)dx]$ is differentiable in $t$ (in a $W^{1,1}$ sense). Next, we can check that the expectation $\mathbb{E}_{\sigma, \mathcal{C}}[v_{i,t}\int_{\mathbb{R}_{+}}u^{\epsilon_{k_{n}}}(t,\,x)f(x)dx]$ converges to $\mathbb{E}_{\sigma, \mathcal{C}}[v_{i,t}\int_{\mathbb{R}_{+}}u^*(t,\,x)f(x)dx]$ for both $i = 1$ and $i = 2$, both weakly and pointwise in $t \in [0, \, T]$, while the limits are also differentiable in $t$ everywhere except the two jump points $t_1$ and $t_2$. This follows because everything is zero outside $[t_1, \, t_2]$, while both $v_{1}$  and $v_{2}$ are constant in $t$ and thus of the form (\ref{d1}) - (\ref{d2}) if we restrict to that interval. Subtracting from each term of (\ref{tested}) the same term but with $u^{\epsilon}$ replaced $u^*$ and then adding it back, we can rewrite this identity as
\begin{eqnarray}\label{testedforlim}
&& \mathbb{E}_{\sigma, \mathcal{C}}\left[V_t\int_{\mathbb{R}_{+}}u^{\epsilon}(t,\,x)f(x)dx\right] \nonumber \\
&& \qquad = \mathbb{E}_{\sigma, \mathcal{C}}\left[V_0\int_{\mathbb{R}_{+}}u_{0}(x)f(x)dx\right] + r\int_{0}^{t}\mathbb{E}_{\sigma, \mathcal{C}}\left[V_s\int_{\mathbb{R}_{+}}u^{\epsilon}(s,\,x)f'(x)dx\right]ds \nonumber \\
&& \qquad \qquad  - \int_{0}^{t}\frac{h^2\left(\sigma_{\frac{s}{\epsilon}}^{1,1}\right)}{2}\Bigg(\mathbb{E}_{\sigma, \mathcal{C}}\left[V_{s}\int_{\mathbb{R}_{+}}u^{\epsilon}(s,\,x)f'(x)dx\right] \nonumber \\
&& \qquad \qquad \qquad \qquad \qquad \qquad \qquad - \mathbb{E}_{\sigma, \mathcal{C}}\left[V_{s}\int_{\mathbb{R}_{+}}u^*(s,\,x)f'(x)dx\right]\Bigg)ds \nonumber \\
&& \qquad \qquad - \int_{0}^{t}\frac{h^2\left(\sigma_{\frac{s}{\epsilon}}^{1,1}\right)}{2}\mathbb{E}_{\sigma, \mathcal{C}}\left[V_{s}\int_{\mathbb{R}_{+}}u^*(s,\,x)f'(x)dx\right]ds \nonumber \\
&& \qquad \qquad + \int_{0}^{t}\frac{h^2\left(\sigma_{\frac{s}{\epsilon}}^{1,1}\right)}{2}\Bigg(\mathbb{E}_{\sigma, \mathcal{C}}\left[V_{s}\int_{\mathbb{R}_{+}}u^{\epsilon}(s,\,x)f''(x)dx\right] \nonumber \\
&& \qquad \qquad \qquad \qquad \qquad \qquad \qquad - \mathbb{E}_{\sigma, \mathcal{C}}\left[V_{s}\int_{\mathbb{R}_{+}}u^*(s,\,x)f''(x)dx\right]\Bigg)ds \nonumber \\
&& \qquad \qquad + \int_{0}^{t}\frac{h^2\left(\sigma_{\frac{s}{\epsilon}}^{1,1}\right)}{2}\mathbb{E}_{\sigma, \mathcal{C}}\left[V_{s}\int_{\mathbb{R}_{+}}u^*(s,\,x)f''(x)dx\right]ds \nonumber \\
&& \qquad \qquad +\int_{0}^{t}\mathbb{E}_{\sigma, \mathcal{C}}\left[v_{1,s}\int_{\mathbb{R}_{+}}
u^{\epsilon}(s,\,x)f(x)dx\right]ds \nonumber \\
&& \qquad \qquad + \rho_{1,1}\int_{0}^{t}h\left(\sigma_{\frac{s}{\epsilon}}^{1,1}\right)\Bigg(\mathbb{E}_{\sigma, \mathcal{C}}\left[v_{2,s}\int_{\mathbb{R}_{+}}u^{\epsilon}(s,\,x)f'(x)dx\right] \nonumber \\
&& \qquad \qquad \qquad \qquad \qquad \qquad \qquad - \mathbb{E}_{\sigma, \mathcal{C}}\left[v_{2,s}\int_{\mathbb{R}_{+}}u^*(s,\,x)f'(x)dx\right]\Bigg)ds \nonumber \\
&& \qquad \qquad + \rho_{1,1}\int_{0}^{t}h\left(\sigma_{\frac{s}{\epsilon}}^{1,1}\right)\mathbb{E}_{\sigma, \mathcal{C}}\left[v_{2,s}\int_{\mathbb{R}_{+}}u^*(s,\,x)f'(x)dx\right]ds. \nonumber \\
\end{eqnarray}
Then we have
\begin{eqnarray}
&& \Bigg|\int_{0}^{t}h\left(\sigma_{\frac{s}{\epsilon}}^{1,1}\right)\Bigg(\mathbb{E}_{\sigma, \mathcal{C}}\left[v_{2,s}\int_{\mathbb{R}_{+}}u^{\epsilon}(s,\,x)f'(x)dx\right] \nonumber \\
&& \qquad \qquad \qquad \qquad \qquad - \mathbb{E}_{\sigma, \mathcal{C}}\left[v_{2,s}\int_{\mathbb{R}_{+}}u^*(s,\,x)f'(x)dx\right]\Bigg)ds \Bigg| \nonumber \\
&& \qquad \leq C\int_{0}^{t}\Bigg|\mathbb{E}_{\sigma, \mathcal{C}}\left[v_{2,s}\int_{\mathbb{R}_{+}}u^{\epsilon}(s,\,x)f'(x)dx\right] \nonumber \\
&& \qquad \qquad \qquad \qquad \qquad  - \mathbb{E}_{\sigma, \mathcal{C}}\left[v_{2,s}\int_{\mathbb{R}_{+}}u^*(s,\,x)f'(x)dx\right]\Bigg|ds, \nonumber
\end{eqnarray}
which tends to zero (when $\epsilon = \epsilon_{k_{n}}$ and $n \rightarrow \infty$) by the dominated convergence theorem, since the quantity inside the last integral converges pointwise to zero and it can be dominated by using (\ref{bound}). The same argument is used to show that the 4th and 6th terms in (\ref{testedforlim}) tend also to zero along the same subsequence. Finally, for any 
term of the form 
\begin{eqnarray}
\int_{0}^{t}h^{p}\left(\sigma_{\frac{s}{\epsilon}}^{1,1}\right)\mathbb{E}_{\sigma, \mathcal{C}}\left[V_{s}\int_{\mathbb{R}_{+}}u^*(s,\,x)f^{(m)}(x)dx\right]ds \nonumber
\end{eqnarray}
for $p, m \in \{0, \, 1, \, 2\}$, we can recall the differentiability of the second factor inside the integral (which was mentioned earlier) and then use integration by parts to write it as:
\begin{eqnarray}
&& \int_{0}^{t}h^{p}\left(\sigma_{\frac{w}{\epsilon}}^{1,1}\right)dw\Bigg(\mathbb{E}_{\sigma, \mathcal{C}}\left[V_{s}\int_{\mathbb{R}_{+}}u^*(t,\,x)f^{(m)}(x)dx\right]\Bigg) \nonumber \\
&& \qquad -\int_{0}^{t}\int_{0}^{s}h^{p}\left(\sigma_{\frac{w}{\epsilon}}^{1,1}\right)dw\Bigg(\mathbb{E}_{\sigma, \mathcal{C}}\left[V_{s}\int_{\mathbb{R}_{+}}u^*(s,\,x)f^{(m)}(x)dx\right]\Bigg)'ds \nonumber
\end{eqnarray}
which converges, by the positive recurrence property, to the quantity
\begin{eqnarray}
&& t\mathbb{E}\left[h^{p}\left(\sigma^{1,1,1,*}\right) \, | \, \mathcal{C}\right]\Bigg( \mathbb{E}_{\sigma, \mathcal{C}}\left[V_{s}\int_{\mathbb{R}_{+}}u^*(t,\,x)f^{(m)}(x)dx\right]\Bigg) \nonumber \\
&& \qquad -\int_{0}^{t}s\mathbb{E}\left[h^{p}\left(\sigma^{1,1,1,*}\right) \, | \, \mathcal{C}\right]\Bigg(\mathbb{E}_{\sigma, \mathcal{C}}\left[V_{s}\int_{\mathbb{R}_{+}}u^*(s,\,x)f^{(m)}(x)dx\right]\Bigg)'ds. \nonumber
\end{eqnarray}
Using integration by parts once more, this last expression is equal to
\begin{eqnarray*}
\mathbb{E}\left[h^{p}\left(\sigma^{1,1,1,*}\right) \, | \, \mathcal{C}\right]\int_{0}^{t}\mathbb{E}_{\sigma, \mathcal{C}}\left[V_{s}\int_{\mathbb{R}_{+}}u^*(s,\,x)f^{(m)}(x)dx\right]ds.
\end{eqnarray*}
This last convergence result holds also if we replace $V$ by $v_{1}$ or $v_{2}$, as we can show by following exactly the same steps in the subinterval $[t_1, \, t_2]$ (where $v_{i}$ is supported for $i \in \{1, \, 2\}$ and where we have differentiability that allows integration by parts). 

If we set now $\epsilon = \epsilon_{k_{n}}$ in (\ref{testedforlim}), take $n \rightarrow +\infty$, and substitute all the above convergence results, we obtain 
\begin{eqnarray}\label{tlspde}
&& \mathbb{E}_{\sigma, \mathcal{C}}\left[V_t\int_{\mathbb{R}_{+}}u^*(t,\,x)f(x)dx\right] \nonumber \\
&& \qquad = \mathbb{E}_{\sigma, \mathcal{C}}\left[V_0\int_{\mathbb{R}_{+}}u_{0}(x)f(x)dx\right] \nonumber \\
&& \qquad \qquad + \left(r-\frac{\sigma_{2,1}^{2}}{2}\right)\int_{0}^{t}\mathbb{E}_{\sigma, \mathcal{C}}\left[V_s\int_{\mathbb{R}_{+}}u^*(s,\,x)f'(x)dx\right]ds \nonumber \\
&& \qquad \qquad + \frac{\sigma_{2,1}^{2}}{2}\int_{0}^{t}
\mathbb{E}_{\sigma, \mathcal{C}}\left[V_s\int_{\mathbb{R}_{+}}u^*(s,\,x)f''(x)dx\right]ds \nonumber \\
&& \qquad \qquad +\int_{0}^{t}\mathbb{E}_{\sigma, \mathcal{C}}\left[v_{1,s}\int_{\mathbb{R}_{+}}
u^*(s,\,x)f(x)dx\right]ds \nonumber \\
&& \qquad \qquad  +\rho_{1,1}\sigma_{1,1}\int_{0}^{t}\mathbb{E}_{\sigma, \mathcal{C}}\left[v_{2,s}\int_{\mathbb{R}_{+}}u^*(s,\,x)f'(x)dx\right]ds.
\end{eqnarray}
Since $\tilde{\mathbb{V}}$ is dense in $\mathbb{V}$, for a fixed $t \leq T$, we can have (\ref{tlspde}) for any square-integrable martingale $\{V_{s}: 0 \leq s \leq t\}$, for which we have $v_{1,s} = 0$ for all $0 \leq s \leq t$. Next, we denote by $R_u(t, \, x)$ the RHS of (\ref{eq:lspde}). Using then Ito's formula for the product of $\int_{\mathbb{R}_{+}}R_u(s,\,x)f(x)dx$ with $V_{s}$ at $s = t$, subtracting $V_t\int_{\mathbb{R}_{+}}u^*(t,\,x)f(x)dx$ from both sides, taking expectations and finally substituting from (\ref{tlspde}), we find that
\begin{eqnarray}
\mathbb{E}_{\sigma, \mathcal{C}}\Bigg[V_t\left(\int_{\mathbb{R}_{+}}R_u(t,\,x)f(x)dx - \int_{\mathbb{R}_{+}}u^*(t,\,x)f(x)dx\right)\Bigg] = 0 \nonumber
\end{eqnarray}
for our fixed $t \leq T$. Using the martingale representation theorem, $V_s$ can be taken equal to $\mathbb{E}_{\sigma, \mathcal{C}}\left[\mathbb{I}_{\mathcal{E}_s} \, | \, \sigma\left(\{W_{s'}^0: s' \leq s\}\right)\right]$ for all $s \leq t$, where we define 
\begin{eqnarray*}
\mathcal{E}_t = \left\{\omega \in \Omega: \int_{\mathbb{R}_{+}}R_u(t,\,x)f(x)dx > \int_{\mathbb{R}_{+}}u^*(t,\,x)f(x)dx\right\},
\end{eqnarray*}
and this implies $V_t = \mathbb{I}_{\mathcal{E}_s}$ allowing us to write
\begin{eqnarray}
\mathbb{E}_{\sigma, \mathcal{C}}\left[\mathbb{I}_{\mathcal{E}_t}\left(\int_{\mathbb{R}_{+}}R_u(t,\,x)f(x)dx - \int_{\mathbb{R}_{+}}u^*(t,\,x)f(x)dx\right)\right] = 0 \nonumber
\end{eqnarray}
for any $0 \leq t \leq T$. If we integrate the above for $t \in [0, \, T]$ we obtain that
\begin{eqnarray}
\int_{0}^{T}\mathbb{E}_{\sigma, \mathcal{C}}\Bigg[\mathbb{I}_{\mathcal{E}_t}\left(\int_{\mathbb{R}_{+}}R_u(t,\,x)f(x)dx - \int_{\mathbb{R}_{+}}u^*(t,\,x)f(x)dx\right)\Bigg]dt = 0 \nonumber
\end{eqnarray}
where the quantity inside the expectation is always non-negative and becomes zero only when $\mathbb{I}_{\mathcal{E}_t} = 0$. This implies $\int_{\mathbb{R}_{+}}R_u(t,\,x)f(x)dx \leq \int_{\mathbb{R}_{+}}u^*(t,\,x)f(x)dx$ almost everywhere, and working in the same way with the indicator of the complement $\mathbb{I}_{\mathcal{E}_t^{c}}$ we can deduce the opposite inequality as well. Thus, we must have $\int_{\mathbb{R}_{+}}R_u(t,\,x)f(x)dx = \int_{\mathbb{R}_{+}}u^*(t,\,x)f(x)dx$ almost everywhere, and since the function $f$ is an arbitrary smooth function with compact support, we can deduce that $R_u$ coincides with $u^*$ almost everywhere, which gives (\ref{eq:lspde}).

If $h$ is bounded from below, we can use (\ref{eq:4.3}) to obtain a uniform (independent of
$\epsilon$) bound for the $H_{0}^{1}(\mathbb{R}_{+})\otimes L_{\sigma, \mathcal{C}}^2(\Omega \times [0, \, T])$ norm of $u^{\epsilon_{k_{n}}}$, which implies that in a further subsequence, the weak convergence to $u^*$ holds also in that Sobolev space, in which (\ref{eq:lspde}) has a unique solution \cite{BHHJR}. This implies convergence of $u^{\epsilon}$ to the unique solution of (\ref{eq:lspde}) in $H_{0}^{1}(\mathbb{R}_{+})\otimes L_{\sigma, \mathcal{C}}^2(\Omega \times [0, \, T])$, as $\epsilon \rightarrow 0^{+}$. The proof is now complete.
\qed
\end{proof}

\begin{proof}[\textbf{Proposition 2.7}]
The upper bound can be obtained by a simple Cauchy-Schwarz inequality,
\begin{eqnarray*}
\tilde{\sigma} &=& \sqrt{\mathbb{E}\Big[h\left(\sigma^{1, 2, 1, *}\right)h\left(\sigma^{1, 2, 2, *}\right)\Big]} \nonumber \\
&\leq& \sqrt{\sqrt{\mathbb{E}\Big[h^2\left(\sigma^{1, 2, 1, *}\right)\Big]}\sqrt{\mathbb{E}\Big[h^2\left(\sigma^{1, 2, 2, *}\right)\Big]}} \nonumber \\
&=& \sqrt{\sigma_{2,1} \times \sigma_{2,1}}. \nonumber \\
&=& \sigma_{2,1}
\end{eqnarray*}
This calculation shows that this bound is only attainable when $\sigma^{i, j, 1, *} = \sigma^{i, j, 2, *}$ for all $i$ and $j$ with $i \neq j$, and this happens only when all the assets share a common stochastic volatility (i.e $\rho_{2} = 1$).

For the lower bound, considering our volatility processes for $i=1$ and $i=2$ started from their $1$-dimensional stationary distributions independently, we have for any $t, \epsilon \geq 0$ 
\begin{eqnarray}\label{lowbound}
&& \mathbb{E}\bigg[\frac{1}{t}\int_{0}^{t}h\Big(\sigma_{\frac{s}{\epsilon}}^{1, 1}\Big)h\left(\sigma_{\frac{s}{\epsilon}}^{2, 1}\right)ds\bigg] \nonumber \\ 
&& \quad = \frac{1}{t}\int_{0}^{t}\mathbb{E}\Big[h\left(\sigma_{\frac{s}{\epsilon}}^{1, 1}\right)h\left(\sigma_{\frac{s}{\epsilon}}^{2, 1}\right)\Big]ds \nonumber \\ 
&& \quad = \frac{1}{t}\int_{0}^{t}\mathbb{E}\Big[h\left(\sigma_{\frac{s}{\epsilon}}^{1, 1}\right)\Big]\mathbb{E}\Big[h\left(\sigma_{\frac{s}{\epsilon}}^{2, 1}\right)\Big]ds \nonumber \\
&& \qquad \quad + \frac{1}{t}\int_{0}^{t}\mathbb{E}\Bigg[\bigg(h\left(\sigma_{\frac{s}{\epsilon}}^{1, 1}\right) - \mathbb{E}\Big[h\left(\sigma_{\frac{s}{\epsilon}}^{1, 1}\right)\Big]\bigg)\bigg(h\left(\sigma_{\frac{s}{\epsilon}}^{2, 1}\right) - \mathbb{E}\Big[h\left(\sigma_{\frac{s}{\epsilon}}^{2, 1}\right)\Big]\bigg)\Bigg]ds \nonumber \\
&& \quad = \sigma_{1,1}^2 \nonumber \\
&& \qquad \quad + \frac{1}{t}\int_{0}^{t}\mathbb{E}\Bigg[\mathbb{E}\bigg[\bigg(h\left(\sigma_{\frac{s}{\epsilon}}^{1, 1}\right) - \mathbb{E}\Big[h\left(\sigma_{\frac{s}{\epsilon}}^{1, 1}\right)\Big]\bigg) \nonumber \\
&& \qquad \qquad \qquad \times \bigg(h\left(\sigma_{\frac{s}{\epsilon}}^{2, 1}\right) - \mathbb{E}\Big[h\left(\sigma_{\frac{s}{\epsilon}}^{1, 1}\right)\Big]\bigg) \, \bigg| \, B^{0}\bigg]\Bigg]ds \nonumber \\
&& \quad = \sigma_{1,1}^2 + \frac{1}{t}\int_{0}^{t}\mathbb{E}\Bigg[\mathbb{E}\bigg[\bigg(h\left(\sigma_{\frac{s}{\epsilon}}^{1, 1}\right) - \sigma_{1,1}\bigg) \, \bigg| \, B^{0}\bigg]^2\Bigg]ds \nonumber \\
&& \quad \geq \sigma_{1,1}^2,
\end{eqnarray}
since $\sigma^{1,1}$ and $\sigma^{2,1}$ are identically distributed, and also independent when $B^{0}$ is given. Taking $\epsilon \rightarrow 0^{+}$ on (\ref{lowbound}) and recalling the positive recurrence property, the definition of $\tilde{\sigma}$, and the dominated convergence theorem on the LHS (since the quantity inside the expectation there is bounded by the square of an upper bound of $h$), we obtain the lower bound, i.e $\tilde{\sigma} \geq \sigma_{1,1}$, which can also be shown to be unattainable in general. Indeed, if we choose $h$ such that its composition $\tilde{h}$ with the square function is strictly increasing and convex, and if $g$ is chosen to be a square root function (thus we are in the CIR volatility case), for any $\alpha > 0$ we have
\begin{eqnarray*}
&& \frac{1}{t}\int_{0}^{t}\mathbb{E}\Bigg[\mathbb{E}\bigg[\bigg(h\Big(\sigma_{\frac{s}{\epsilon}}^{1, 1}\Big) - \sigma_{1,1}\bigg) \, \bigg| \, B^{0}\bigg]^2\Bigg]ds \nonumber \\
&& \qquad = \mathbb{E}\Bigg[\frac{1}{t}\int_{0}^{t}\Bigg(\mathbb{E}\left[\tilde{h}\bigg(\sqrt{\sigma_{\frac{s}{\epsilon}}^{1, 1}}\right)\, \bigg| \, B^{0}\bigg]- \sigma_{1,1}\Bigg)^2ds\Bigg] \nonumber \\
&& \qquad \geq \alpha^2\mathbb{E}\left[\frac{1}{t}\int_{0}^{t}\mathbb{I}_{\sigma_{\frac{s}{\epsilon}}^{B^{0}, h} \geq \alpha + \sigma_{1,1}}ds\right]
\end{eqnarray*}
where $\sigma_{s}^{B^{0}, h} := \mathbb{E}[\tilde{h}(\sqrt{\sigma_{s}^{1, 1}})\, | \, B^{0}] \geq \tilde{h}(\sigma_{s}^{B^{0}})$ for $\sigma_{s}^{B^{0}} := \mathbb{E}[\sqrt{\sigma_{s}^{1, 1}}\, | \, B^{0}]$, which implies that
\begin{eqnarray}\label{erg2}
&& \frac{1}{t}\int_{0}^{t}\mathbb{E}\Bigg[\mathbb{E}\bigg[\bigg(h\left(\sigma_{\frac{s}{\epsilon}}^{1, 1}\right) - \sigma_{1,1}\bigg) \, \bigg| \, B^{0}\bigg]^2\Bigg]ds \nonumber \\
&& \qquad \geq \alpha^2\mathbb{E}\left[\frac{1}{t}\int_{0}^{t}\mathbb{I}_{\sigma_{\frac{s}{\epsilon}}^{B^{0}} \geq \tilde{h}^{-1}\left(\alpha + \sigma_{1,1}\right)}ds\right].
\end{eqnarray}
Let $\sigma_{t}^{\rho}$ be the solution to the SDE
\begin{eqnarray}
\sigma_{t}^{\rho} = \sigma_{0}^{B^{0}} + \frac{1}{2}\int_{0}^{t}\left(\kappa\theta - \frac{v^2}{4}\right)\frac{1}{\sigma_{s}^{\rho}}ds + \frac{\kappa}{2}\int_{0}^{t}\sigma_{s}^{\rho}ds + \frac{\rho_{2}v}{2}B_{s}^{0}. \nonumber
\end{eqnarray}
Then $\sigma^{\rho}$ can be shown to be the square root of a CIR process having the same mean-reversion and vol-of-vol as $\sigma^{1,1}$ and a different stationary mean, which satisfies the Feller condition for not hitting zero at a finite time. If for some $t_1 > 0$ we have $\sigma_{t_1}^{\rho} > \sigma_{t_1}^{B^{0}}$, we consider $t_0 = \sup\{s \leq t_1: \sigma_{s}^{\rho} = \sigma_{s}^{B^{0}}\}$ which is obviously non-negative. Then, since $\mathbb{E}[{\frac{1}{\sqrt{\sigma_{s}^{1,1}}}}\, | \, B^{0}] \geq \frac{1}{\mathbb{E}[\sqrt{\sigma_{s}^{1,1}}\, | \, B^{0}]} = \frac{1}{\sigma_{s}^{B^{0}}}$ we have 
\begin{eqnarray}
\sigma_{t_1}^{B^{0}} &=& \sigma_{t_0}^{B^{0}} + \frac{1}{2}\int_{t_0}^{t_1}\left(\kappa\theta - \frac{v^2}{4}\right)\mathbb{E}\Bigg[{\frac{1}{\sqrt{\sigma_{s}^{1,1}}}}\, \Bigg| \, B^{0}\Bigg]ds \nonumber \\ 
&& \qquad - \frac{\kappa}{2}\int_{t_0}^{t_1}\sigma_{s}^{B^{0}}ds + \frac{\rho_{2}v}{2}\left(B_{t_1}^{0} - B_{t_0}^{0}\right) \nonumber \\
&\geq& \sigma_{t_0}^{B^{0}} + \frac{1}{2}\int_{t_0}^{t_1}\left(\kappa\theta - \frac{v^2}{4}\right)\frac{1}{\sigma_{s}^{B^{0}}}ds - \frac{\kappa}{2}\int_{t_0}^{t_1}\sigma_{s}^{B^{0}}ds + \frac{\rho_{2}v}{2}\left(B_{t_1}^{0} - B_{t_0}^{0}\right) \nonumber \\
&\geq& \sigma_{t_0}^{\rho} + \frac{1}{2}\int_{t_0}^{t_1}\left(\kappa\theta - \frac{v^2}{4}\right)\frac{1}{\sigma_{s}^{\rho}}ds - \frac{\kappa}{2}\int_{t_0}^{t_1}\sigma_{s}^{\rho}ds + \frac{\rho_{2}v}{2}\left(B_{t_1}^{0} - B_{t_0}^{0}\right) \nonumber \\
&=& \sigma_{t_1}^{\rho} \nonumber 
\end{eqnarray}
which is a contradiction. Thus $\sigma_{s}^{\rho} \leq \sigma_{s}^{B^{0}}$ for all $s \geq 0$, and in (\ref{erg2}) this gives 
\begin{eqnarray*}
&& \frac{1}{t}\int_{0}^{t}\mathbb{E}\Bigg[\mathbb{E}\bigg[\bigg(h\Big(\sigma_{\frac{s}{\epsilon}}^{1, 1}\Big) - \sigma_{1,1}\bigg) \, \bigg| \, B_{\cdot}^{0}\bigg]^2\Bigg]ds \nonumber \\
&& \qquad \geq \alpha^2\mathbb{E}\left[\frac{1}{t}\int_{0}^{t}\mathbb{I}_{\sigma_{\frac{s}{\epsilon}}^{\rho} \geq \tilde{h}^{-1}\left(\alpha + \sigma_{1,1}\right)}ds\right]. 
\end{eqnarray*}
By the positive recurrence of $\sigma^{\rho}$ (which is the root of a CIR process, the ergodicity of which has been discussed in \cite{FOU}), the RHS of the above converges to $\alpha^2\mathbb{P}(\sigma^{\rho, *} \geq \tilde{h}^{-1}(\alpha + \sigma_{1,1}))$ as $\epsilon \rightarrow 0^{+}$, where $\sigma^{\rho, *}$ has the stationary distribution of $\sigma^{\rho}$. This expression can only be zero when $\sigma^{\rho, *}$ is a constant, and since the square of $\sigma^{\rho}$ satisfies Feller's boundary condition, this can only happen when $\rho_{2} = 0$. In that case, we can easily check that $\sigma^{1, 2, 1, *}$ and $\sigma^{1, 2, 1, *}$ are independent, which implies that $\tilde{\sigma} = \sigma_{1,1}$. This completes the proof.
\qed
\end{proof}

\begin{proof}[\textbf{Corollary 2.8}] 
Suppose that $\mathbb{P}(X_t^{1, \epsilon} \in \mathcal{I} \, | \, W_{\cdot}^{0}, \, B_{\cdot}^{0}, \, \mathcal{G})$ converges to $\mathbb{P}(X_t^{1, *} \in \mathcal{I} \, | \, W_{\cdot}^{0}, \, \mathcal{G})$ in probability, under the assumptions of both Theorem~\ref{theorem} and Theorem~\ref{thm1}. The same convergence 
has to hold in a strong $L^2$ sense for some sequence $\epsilon_n \downarrow 0$, since it
will hold $\mathbb{P}$ - almost surely for some sequence, and then we can apply the dominated convergence theorem. Therefore, the same convergence must hold weakly in $L^2$ as well. However, assuming for simplicity that $(r_i, \, \rho_{1,i})$ is also a constant vector $(r, \, \rho_{1})$ for all $i$ and fixing a sufficiently integrable and $\sigma(W_{\cdot}^{0}, \, B_{\cdot}^{0}) \cap \mathcal{G}$-measurable random variable $\Xi$, by Theorem~\ref{thm1} we have
\begin{eqnarray*}
&&\lim_{n \rightarrow +\infty}\mathbb{E}\bigg[\Xi\mathbb{P}\Big[X_t^{1, \epsilon_n} \in \mathcal{I} \, \Big| \, W^{0}, \, B^{0}, \, \mathcal{G}\Big] \bigg] \nonumber \\
&& \qquad = \lim_{n \rightarrow +\infty}\mathbb{E}\bigg[\Xi\mathbb{P}\Big[X_t^{1, \epsilon_n} \in \mathcal{I} \, \Big| \, W^{0}, \,\sigma^{1,1}, \, \mathcal{G}\Big] \bigg] \nonumber \\
&& \qquad = \lim_{n \rightarrow +\infty}\mathbb{E}\left[\int_{0}^{+\infty}\Xi\mathbb{I}_{\mathcal{I}}(x)u^{\epsilon_n}(t, x)dx\right] \nonumber \\
&& \qquad = \mathbb{E}\left[\int_{0}^{+\infty}\Xi\mathbb{I}_{\mathcal{I}}(x)u^*(t, x)dx\right] \nonumber \\
&& \qquad = \mathbb{E}\bigg[\Xi\mathbb{P}\Big[X_t^{1, w} \in \mathcal{I} \, \Big| \, W^{0}, \, \mathcal{G}\Big]\bigg], \nonumber \\
\end{eqnarray*}
where for each $i$ we define 
\[
\begin{array}{rcl}
X_{t}^{i, w} &=& x^{i} + \left(r - \frac{\sigma_{2,1}^2}{2}\right)t + \rho_{1}'\sigma_{2,1}W_{t}^{0} + \sqrt{1 - \left(\rho_{1}'\right)^2}\sigma_{2,1}W_{t}^{i}, \,\, 0 \leq t \leq T_{i}^w \\
X_{t}^{i, w} &=& 0, \, t \geq T_{i}^w\\
T_{i}^w &=& \inf\{ t\geq 0: X_t^{i, w}=0\}, 
\end{array}
\]
with $\rho_{1}' = \rho_{1}\frac{\sigma_{1,1}}{\sigma_{2,1}}$, in which the density of $X_t^{i, w}$ given $W^0$ and $\mathcal{G}$ is the unique solution $u^*$ to (\ref{eq:lspde}) \cite{HL16}. Therefore, by the uniqueness of a weak limit we must have $\mathbb{P}[X_t^{1, *} \in \mathcal{I} \, | \, W^{0}, \, \mathcal{G}] = \mathbb{P}[X_t^{1, w} \in \mathcal{I} \, | \, W^{0}, \, \mathcal{G}]$ $\mathbb{P}$-almost surely, which cannot be true for any interval $\mathcal{I}$, as otherwise the processes $X_{\cdot}^{1, w}$ and $X_{\cdot}^{1, *}$ would coincide, which is clearly not the case here. Indeed, this can only be true when $\tilde{\rho}_{1,1} = \rho_1' \Leftrightarrow \tilde{\sigma} = \sigma_{1,1}$, and by Proposition~\ref{prp3} this is generally not the case unless $\rho_2 = 0$. 
\qed
\end{proof}

\section{Proofs: small vol-of-vol setting}

We proceed now to the proofs of Proposition~\ref{lem1}, Theorem~\ref{thm5.2} and Corollary~\ref{rateofcon}, the main results of Section 3.

\begin{proof}[\bf{Proposition 3.1}]
First, we will show that each volatility process has a finite $2p$-moment for any $p \in \mathbb{N}$. Indeed, we fix a $p \in \mathbb{N}$ and we consider the sequence of stopping times $\{\tau_{n, \epsilon}: n \in \mathbb{N}\}$, where $\tau_{n, \epsilon} = \inf\{t \geq 0: \sigma_{t}^{i, \epsilon} > n\}$. Setting $\sigma_{t}^{i, n, \epsilon} = \sigma_{t \wedge \tau_{n, \epsilon}}^{i, \epsilon}$, by Ito's formula we have
\begin{eqnarray}\label{convtomean}
\left(\sigma_{t}^{i, n, \epsilon} - \theta_i\right)^{2p} &=& \left(\sigma_{0}^{i, n, \epsilon} - \theta_i\right)^{2p} - \frac{2p\kappa_i}{\epsilon}\int_{0}^{t}\mathbb{I}_{\left[0, \tau_{n, \epsilon}\right]}(s)\left(\sigma_{s}^{i, n, \epsilon} - \theta_i\right)^{2p}ds \nonumber \\
&& + 2p\xi_i\int_{0}^{t}\mathbb{I}_{\left[0, \tau_{n, \epsilon}\right]}(s)\left(\sigma_{s}^{i, n, \epsilon} - \theta_i\right)^{2p-1}g\left(\sigma_{s}^{i, n, \epsilon}\right)d\tilde{B}_{s}^{i} \nonumber \\
&& + p(2p-1)\xi_i^2\int_{0}^{t}\mathbb{I}_{\left[0, \tau_{n, \epsilon}\right]}(s)\left(\sigma_{s}^{i, n, \epsilon} - \theta_i\right)^{2p-2}g^2\left(\sigma_{s}^{i, n, \epsilon}\right)ds \nonumber \\
\end{eqnarray}
for $\tilde{B}_{s}^{i} = \sqrt{1-\rho_{2,i}^{2}}B_{t}^{i}+\rho_{2,i}B_{t}^{0}$, where the stochastic integral is a martingale. Taking expectations, setting $f(t, n, p, \epsilon) = \mathbb{E}[(\sigma_{t}^{i, n, \epsilon} - \theta_i)^{2p}]$ and using the growth condition of $g$ ($|g(x)| \leq C_{1,g} + C_{2,g}|x|$ for all $x \in \mathbb{R}$) and simple inequalities, we can easily obtain
\begin{eqnarray}
f(t, n, p, \epsilon) \leq M + M'\int_{0}^{t}f(s, n, p, \epsilon)ds \nonumber
\end{eqnarray} 
with $M, M'$ depending only on $p$, $c_g$ and the bounds of $\sigma^{i}, \xi_i, \theta_i$. Thus, using Gronwall's inequality we get a uniform (in $n$) estimate for $f(t, n, p, \epsilon)$, and then by Fatou's lemma we obtain the desired finiteness of
\begin{eqnarray*}
f(t, p, \epsilon) := \mathbb{E}\left[\left(\sigma_{t}^{i, \epsilon} - \theta_i\right)^{2p}\right].
\end{eqnarray*}
This implies the almost sure finiteness of the conditional expectation
\begin{eqnarray*}
f_{\mathcal{C}}(t, p, \epsilon) := \mathbb{E}\bigg[\left(\sigma_{t}^{i, \epsilon} - \theta_i\right)^{2p} \, \bigg| \, \mathcal{C}\bigg]
\end{eqnarray*}
as well. 

Taking expectations given $\mathcal{C}$, letting $n \rightarrow +\infty$ on (\ref{convtomean}), using the monotone convergence theorem (all quantities are monotone for large enough $n$) and the growth condition on $g$, we find that
\begin{eqnarray}
f_{\mathcal{C}}(t, p, \epsilon) \leq M + \left(M' - \frac{2\kappa_i}{\epsilon}\right)\int_{0}^{t}f_{\mathcal{C}}(s, p, \epsilon)ds, \nonumber
\end{eqnarray} 
where again, $M, M'$ depend only on $p$, $c_g$ and the bounds of $\sigma^{i}, \xi_i, \theta_i$. Using Grownwall's inequality again on the above, we obtain the estimate
\begin{eqnarray}
\int_{0}^{t}f_{\mathcal{C}}(s, p, \epsilon)ds &\leq& M\int_{0}^{t}e^{\left(M' - \frac{2\kappa_i}{\epsilon}\right)(t-s)}ds \nonumber \\ &\leq& Me^{M't}\int_{0}^{t}e^{\left(- \frac{2\kappa_i}{\epsilon}\right)(t-s)}ds \nonumber \\
&=& \frac{\epsilon}{2\kappa_i}e^{M't}\left(1 - e^{\left(- \frac{2\kappa_i}{\epsilon}\right)t}\right) \nonumber \\
&<& \frac{\epsilon}{2c_{\kappa}}e^{M't}. \nonumber
\end{eqnarray}
Then, we have that $\int_{0}^{t}f(s, p, \epsilon)ds = \mathbb{E}\left[\int_{0}^{t}f_{\mathcal{C}}(s, p, \epsilon)ds\right] < \frac{\epsilon}{2c_{\kappa}}e^{M't}$, and this gives the desired result.
\qed
\end{proof}

\begin{proof}[\textbf{Theorem 3.2}]
We can easily check that $v^{0, \epsilon}$ and $v^{0}$ are the unique solutions to the SPDEs (\ref{eq:5.2}) and (\ref{eq:5.2lim}) respectively in $L^2(\Omega \times [0, \, T] ; H^{2}(\mathbb{R}_{+}))$, under the boundary conditions $v_{x}^{0, \epsilon}(t, \, 0) = 0$ and $v_{x}^{0}(t, \, 0) = 0$ respectively. Subtracting the SPDEs satisfied by $v^{0, \epsilon}$ and $v^{0}$ and setting $v^{d, \epsilon} = v^{0} - v^{0, \epsilon}$, we can easily verify that
\begin{eqnarray}
&&v^{d,\epsilon}\left(t, \, x\right) \nonumber \\
&& \qquad = -\frac{1}{2}\int_{0}^{t}\Big(h^2\left(\sigma_{s}^{1, \epsilon}\right) - h^2\left(\theta_1\right) \Big)v_{x}^{0, \epsilon}\left(s, \, x\right)ds \nonumber \\
&& \qquad \qquad + \int_{0}^{t}\left(r - \frac{h^2\left(\theta_1\right)}{2} \right)v_{x}^{d, \epsilon}\left(s, \, x\right)ds \nonumber \\
&& \qquad\qquad + \frac{1}{2}\int_{0}^{t}\Big(h^2\left(\sigma_{s}^{1, \epsilon}\right) - h^2\left(\theta_1\right) \Big)v_{xx}^{0, \epsilon}\left(s, \, x\right)ds + \int_{0}^{t}\frac{h^2\left(\theta_1\right)}{2}v_{xx}^{d, \epsilon}\left(s, \, x\right)ds \nonumber \\
&& \qquad\qquad + \rho_{1,1}\int_{0}^{t}\Big(h\left(\sigma_{s}^{1, \epsilon}\right) - h\left(\theta_1\right) \Big)v_{x}^{0, \epsilon}\left(s, \, x\right)dW_{s}^{0} \nonumber \\
&& \qquad\qquad + \rho_{1,1}\int_{0}^{t}h\left(\theta_1\right)v_{x}^{d, \epsilon}\left(s, \, x\right)dW_{s}^{0}. \nonumber
\end{eqnarray}
Now using Ito's formula for the $L^2$ norm (see Krylov and Rozovskii \cite{KR81}), given the volatility path 
and $\mathcal{C}$, we obtain
\begin{eqnarray}\label{basicest5}
&&\mathbb{E}_{\sigma, \mathcal{C}}\left[\int_{\mathbb{R}_{+}}\Big(v^{d, \epsilon}\left(t, \, x\right)\Big)^2dx\right] \nonumber \\
&& \qquad \qquad = -\int_{0}^{t}\Big(h^2\left(\sigma_{s}^{1, \epsilon}\right) - h^2\left(\theta_1\right) \Big)\mathbb{E}_{\sigma, \mathcal{C}}\left[\int_{\mathbb{R}_{+}}v_{x}^{0, \epsilon}\left(s, \, x\right)v^{d, \epsilon}\left(t, \, x\right)dx\right]ds \nonumber \\
&& \qquad\qquad\qquad + 2\left(r - \frac{h^2\left(\theta_1\right)}{2} \right)\int_{0}^{t}\mathbb{E}_{\sigma, \mathcal{C}}\left[\int_{\mathbb{R}_{+}}v_{x}^{d, \epsilon}\left(s, \, x\right)v^{d, \epsilon}\left(t, \, x\right)dx\right]ds \nonumber \\
&& \qquad\qquad\qquad - \int_{0}^{t}\Big(h^2\left(\sigma_{s}^{1, \epsilon}\right) - h^2\left(\theta_1\right) \Big) \nonumber \\ 
&& \qquad\qquad\qquad\qquad \times \mathbb{E}_{\sigma, \mathcal{C}}\left[\int_{\mathbb{R}_{+}}v_{x}^{0, \epsilon}\left(s, \, x\right)v_{x}^{d, \epsilon}\left(t, \, x\right)dx\right]ds \nonumber \\
&&\qquad\qquad\qquad - \int_{0}^{t}h^2\left(\theta_1\right)\mathbb{E}_{\sigma, \mathcal{C}}\left[\int_{\mathbb{R}_{+}}\Big(v_{x}^{d, \epsilon}\left(s, \, x\right)\Big)^2dx\right]ds \nonumber \\
&&\qquad\qquad\qquad + \rho_{1,1}^2\int_{0}^{t}\Big(h\left(\sigma_{s}^{1, \epsilon}\right) - h\left(\theta_1\right) \Big)^2\mathbb{E}_{\sigma, \mathcal{C}}\left[\int_{\mathbb{R}_{+}}\Big(v_{x}^{0, \epsilon}\left(s, \, x\right)\Big)^2dx\right]ds \nonumber \\
&&\qquad\qquad\qquad + 2\rho_{1,1}^2h\left(\theta_1\right)\int_{0}^{t}\Big(h\left(\sigma_{s}^{1, \epsilon}\right) - h\left(\theta_1\right) \Big) \nonumber \\
&& \qquad\qquad\qquad\qquad \times \mathbb{E}_{\sigma, \mathcal{C}}\left[\int_{\mathbb{R}_{+}}v_{x}^{0, \epsilon}\left(s, \, x\right)v_{x}^{d, \epsilon}\left(t, \, x\right)dx\right]ds \nonumber \\
&& \qquad\qquad\qquad+ \rho_{1,1}^2\int_{0}^{t}h^2\left(\theta_1\right)\mathbb{E}_{\sigma, \mathcal{C}}\left[\int_{\mathbb{R}_{+}}\Big(v_{x}^{d, \epsilon}\left(s, \, x\right)\Big)^2dx\right]ds + N(t, \epsilon)
\end{eqnarray}
where $N(t, \epsilon)$ is some noise due to the correlation between $B^{0}$ and $W^{0}$, with $\mathbb{E}[N(t, \epsilon)] = 0$. In particular, since for some Brownian motion $V^0$ independent from $B^0$ we could have written $W^0 = \sqrt{1 - \rho_3^2}V^0 + \rho_3B^0$, we will have 
\[N(t, \epsilon) = 2\rho_{1,1}\int_{0}^{t}\mathbb{E}_{\sigma, \, \mathcal{C}}\left[h\left(\sigma_s^{1, \epsilon}\right)\int_{\mathbb{R}_{+}}v^{d, \epsilon}\left(s, \, x\right)v_x^{0, \epsilon}\left(s, \, x\right)dx\right]dB_s^0.\] 
Next, we can apply 2. of Theorem 4.1 in \cite{HK17} to the SPDE (\ref{eq:5.2}) to find  
\[\left\Vert v_{x}^{0, \epsilon}(s, \cdot)\right\Vert_{L_{\sigma, \mathcal{C}}^2\left(\Omega \times \mathbb{R}_{+}\right)} = \left\Vert u^{\epsilon}(s, \cdot)\right\Vert_{L_{\sigma, \mathcal{C}}^2\left(\Omega \times \mathbb{R}_{+}\right)} \leq \left\Vert u_{0}(\cdot)\right\Vert_{L^2\left(\Omega \times \mathbb{R}_{+}\right)} \]
for all $s \geq 0$. Using this expression, we can obtain the following estimate 
\begin{eqnarray}\label{smallestimate1}
&&\int_{0}^{t}\Big(h^2\left(\sigma_{s}^{1, \epsilon}\right) - h^2\left(\theta_1\right) \Big)\mathbb{E}_{\sigma, \mathcal{C}}\left[\int_{\mathbb{R}_{+}}v_{x}^{0, \epsilon}\left(s, \, x\right)v^{d, \epsilon}\left(t, \, x\right)dx\right]ds \nonumber \\
&& \qquad \leq \int_{0}^{t}\Big(h^2\left(\sigma_{s}^{1, \epsilon}\right) - h^2\left(\theta_1\right) \Big)\left\Vert v_{x}^{0, \epsilon}(s, \cdot)\right\Vert_{L_{\sigma, \mathcal{C}}^2\left(\Omega \times \mathbb{R}_{+}\right)} \nonumber \\
&& \qquad \qquad \qquad \qquad \qquad \times \left\Vert v^{d, \epsilon}(s, \cdot)\right\Vert_{L_{\sigma, \mathcal{C}}^2\left(\Omega \times \mathbb{R}_{+}\right)}ds \nonumber \\
&& \qquad \leq \left\Vert u_{0}(\cdot)\right\Vert_{L^2\left(\Omega \times \mathbb{R}_{+}\right)}\sqrt{\int_{0}^{t}\bigg(h^2\left(\sigma_{s}^{1, \epsilon}\right) - h^2\left(\theta_1\right) \bigg)^2ds} \nonumber \\
&& \qquad \qquad \qquad \qquad \qquad \times \sqrt{\int_{0}^{t}\left\Vert v^{d, \epsilon}(s, \cdot)\right\Vert_{L_{\sigma, \mathcal{C}}^2\left(\Omega \times \mathbb{R}_{+}\right)}^2ds} \nonumber \\
&& \qquad \leq \frac{1}{2}\left\Vert u_{0}(\cdot)\right\Vert_{L^2\left(\Omega \times \mathbb{R}_{+}\right)}^2\int_{0}^{t}\Big(h^2\left(\sigma_{s}^{1, \epsilon}\right) - h^2\left(\theta_1\right) \Big)^2ds \nonumber \\
&& \qquad \qquad \qquad + \frac{1}{2}\int_{0}^{t}\left\Vert v^{d, \epsilon}(s, \cdot)\right\Vert_{L_{\sigma, \mathcal{C}}^2\left(\Omega \times \mathbb{R}_{+}\right)}^2ds,
\end{eqnarray}
and in the same way
\begin{eqnarray}\label{smallestimate2}
&&\int_{0}^{t}\Big(h^2\left(\sigma_{s}^{1, \epsilon}\right) - h^2\left(\theta_1\right) \Big)\mathbb{E}_{\sigma, \mathcal{C}}\left[\int_{\mathbb{R}_{+}}v_{x}^{0, \epsilon}\left(s, \, x\right)v_{x}^{d, \epsilon}\left(t, \, x\right)dx\right]ds \nonumber \\
&& \qquad \leq \left\Vert u_{0}(\cdot)\right\Vert_{L^2\left(\Omega \times \mathbb{R}_{+}\right)}\sqrt{\int_{0}^{t}\bigg(h^2\left(\sigma_{s}^{1, \epsilon}\right) - h^2\left(\theta_1\right) \bigg)^2ds}\nonumber \\
&& \qquad \qquad \qquad \qquad \qquad \times \sqrt{\int_{0}^{t}\left\Vert v_{x}^{d, \epsilon}(s, \cdot)\right\Vert_{L_{\sigma, \mathcal{C}}^2\left(\Omega \times \mathbb{R}_{+}\right)}^2ds} \nonumber \\
&& \qquad \leq \frac{1}{2\eta}\left\Vert u_{0}(\cdot)\right\Vert_{L^2\left(\Omega \times \mathbb{R}_{+}\right)}^2\int_{0}^{t}\Big(h^2\left(\sigma_{s}^{1, \epsilon}\right) - h^2\left(\theta_1\right) \Big)^2ds \nonumber \\
&& \qquad \qquad \qquad + \frac{\eta}{2}\int_{0}^{t}\left\Vert v_{x}^{d, \epsilon}(s, \cdot)\right\Vert_{L_{\sigma, \mathcal{C}}^2\left(\Omega \times \mathbb{R}_{+}\right)}^2ds
\end{eqnarray}
and
\begin{eqnarray}\label{smallestimate3}
&&\int_{0}^{t}\Big(h\left(\sigma_{s}^{1, \epsilon}\right) - h\left(\theta_1\right) \Big)\mathbb{E}_{\sigma, \mathcal{C}}\left[\int_{\mathbb{R}_{+}}v_{x}^{0, \epsilon}\left(s, \, x\right)v_{x}^{d, \epsilon}\left(t, \, x\right)dx\right]ds \nonumber \\
&& \qquad \leq \left\Vert u_{0}(\cdot)\right\Vert_{L^2\left(\Omega \times \mathbb{R}_{+}\right)}\sqrt{\int_{0}^{t}\bigg(h\left(\sigma_{s}^{1, \epsilon}\right) - h\left(\theta_1\right) \bigg)^2ds}\nonumber \\
&& \qquad \qquad \qquad \qquad \qquad \times \sqrt{\int_{0}^{t}\left\Vert v_{x}^{d, \epsilon}(s, \cdot)\right\Vert_{L_{\sigma, \mathcal{C}}^2\left(\Omega \times \mathbb{R}_{+}\right)}^2ds} \nonumber \\
&& \qquad \leq \frac{1}{2\eta}\left\Vert u_{0}(\cdot)\right\Vert_{L^2\left(\Omega \times \mathbb{R}_{+}\right)}^{2}\int_{0}^{t}\Big(h\left(\sigma_{s}^{1, \epsilon}\right) - h\left(\theta_1\right) \Big)^2ds \nonumber \\
&& \qquad \qquad \qquad + \frac{\eta}{2}\int_{0}^{t}\left\Vert v_{x}^{d, \epsilon}(s, \cdot)\right\Vert_{L_{\sigma, \mathcal{C}}^2\left(\Omega \times \mathbb{R}_{+}\right)}^2ds
\end{eqnarray}
for some $\eta > 0$. Moreover, we have the estimate
\begin{eqnarray}\label{smallestimate4}
&& \int_{0}^{t}\mathbb{E}_{\sigma, \mathcal{C}}\left[\int_{\mathbb{R}_{+}}v_{x}^{d, \epsilon}\left(s, \, x\right)v^{d, \epsilon}\left(t, \, x\right)dx\right]ds \nonumber \\
&& \qquad \leq \int_{0}^{t}\left\Vert v_{x}^{d, \epsilon}(s, \cdot)\right\Vert_{L_{\sigma, \mathcal{C}}^2\left(\Omega \times \mathbb{R}_{+}\right)}\left\Vert v^{d, \epsilon}(s, \cdot)\right\Vert_{L_{\sigma, \mathcal{C}}^2\left(\Omega \times \mathbb{R}_{+}\right)}ds \nonumber \\
&& \qquad \leq \left(\int_{0}^{t}\left\Vert v^{d, \epsilon}(s, \cdot)\right\Vert_{L_{\sigma, \mathcal{C}}^2\left(\Omega \times \mathbb{R}_{+}\right)}^2ds\right)^{1/2}\left(\int_{0}^{t}\left\Vert v_{x}^{d, \epsilon}(s, \cdot)\right\Vert_{L_{\sigma, \mathcal{C}}^2\left(\Omega \times \mathbb{R}_{+}\right)}^2ds\right)^{1/2} \nonumber \\
&& \qquad \leq \frac{1}{2\eta}\int_{0}^{t}\left\Vert v^{d, \epsilon}(s, \cdot)\right\Vert_{L_{\sigma, \mathcal{C}}^2\left(\Omega \times \mathbb{R}_{+}\right)}^2ds + \frac{\eta}{2}\int_{0}^{t}\left\Vert v_{x}^{d, \epsilon}(s, \cdot)\right\Vert_{L_{\sigma, \mathcal{C}}^2\left(\Omega \times \mathbb{R}_{+}\right)}^2ds \nonumber \\
\end{eqnarray}
and by using $\Vert v_{x}^{0, \epsilon}(s, \cdot)\Vert_{L_{\sigma, \mathcal{C}}^2\left(\Omega \times \mathbb{R}_{+}\right)} \leq \Vert u_{0}(\cdot)\Vert_{L^2(\Omega \times \mathbb{R}_{+})}$ again, we also obtain
\begin{eqnarray}\label{smallestimate5}
&& \int_{0}^{t}\Big(h\left(\sigma_{s}^{1, \epsilon}\right) - h\left(\theta_1\right) \Big)^2\mathbb{E}_{\sigma, \mathcal{C}}\left[\int_{\mathbb{R}_{+}}\Big(v_{x}^{0, \epsilon}\left(s, \, x\right)\Big)^2dx\right]ds \nonumber \\
&& \qquad \leq \left\Vert u_{0}(\cdot)\right\Vert_{L^2\left(\Omega \times \mathbb{R}_{+}\right)}^2\int_{0}^{t}\Big(h\left(\sigma_{s}^{1, \epsilon}\right) - h\left(\theta_1\right) \Big)^2ds.
\end{eqnarray}
Using (\ref{smallestimate1}), (\ref{smallestimate2}), (\ref{smallestimate3}), (\ref{smallestimate4}) and (\ref{smallestimate5}) in (\ref{basicest5}), and then taking $\eta$ to be sufficiently small, we get the estimate
\begin{eqnarray}\label{almostlastestimate}
&&\left\Vert v^{d, \epsilon}(t, \cdot)\right\Vert_{L_{\sigma, \mathcal{C}}^2\left(\Omega \times \mathbb{R}_{+}\right)}^2 + m\int_{0}^{t}\left\Vert v_{x}^{d, \epsilon}(s, \cdot)\right\Vert_{L_{\sigma, \mathcal{C}}^2\left(\Omega \times \mathbb{R}_{+}\right)}^2ds \nonumber \\
&& \qquad \leq M\int_{0}^{t}\left\Vert v^{d, \epsilon}(s, \cdot)\right\Vert_{L_{\sigma, \mathcal{C}}^2\left(\Omega \times \mathbb{R}_{+}\right)}^2ds + N(t, \epsilon) + MH(\epsilon)
\end{eqnarray}
for all $t \in [0, \, T]$, where 
\[H(\epsilon) = \int_{0}^{T}\Big(h^2\left(\sigma_{s}^{1, \epsilon}\right) - h^2\left(\theta_1\right) \Big)^2ds + \int_{0}^{T}\Big(h\left(\sigma_{s}^{1, \epsilon}\right) - h\left(\theta_1\right) \Big)^2ds\]
and where $M, m > 0$ are constants independent of the fixed volatility path. Taking expectations in (\ref{almostlastestimate}) to average over all volatility paths, 
we find that
\begin{eqnarray*}
&&\left\Vert v^{d, \epsilon}(t, \cdot)\right\Vert_{L^2\left(\Omega \times \mathbb{R}_{+}\right)}^2 + m\int_{0}^{t}\left\Vert v_{x}^{d, \epsilon}(s, \cdot)\right\Vert_{L^2\left(\Omega \times \mathbb{R}_{+}\right)}^2ds \nonumber \\
&& \qquad \leq M\int_{0}^{t}\left\Vert v^{d, \epsilon}(s, \cdot)\right\Vert_{L^2\left(\Omega \times \mathbb{R}_{+}\right)}^2ds + M\mathbb{E}\left[H(\epsilon)\right],
\end{eqnarray*}
and using Gronwall's inequality on the above we finally obtain
\begin{eqnarray}\label{lastestimate}
&&\left\Vert v^{d, \epsilon}(t, \cdot)\right\Vert_{L^2\left(\Omega \times \mathbb{R}_{+}\right)}^2 + m\int_{0}^{t}\left\Vert v_{x}^{d, \epsilon}(s, \cdot)\right\Vert_{L^2\left(\Omega \times \mathbb{R}_{+}\right)}^2ds \leq M'\mathbb{E}\left[H(\epsilon)\right] \nonumber 
\end{eqnarray}
for some $M' > 0$, with $\mathbb{E}[H(\epsilon)] = \mathcal{O}(\epsilon)$ as $\epsilon \rightarrow 0^{+}$. This last result follows since for $\tilde{h} \in \{h, \, h^2\}$ we can use the mean value theorem to find that
\begin{eqnarray}\label{polybound}
\int_{0}^{T}\left(\tilde{h}\left(\sigma_{s}^{1, \epsilon}\right) - \tilde{h}\left(\theta_1\right) \right)^2ds &=& \int_{0}^{T}\tilde{h}'\left(\sigma_{s,*}^{1, \epsilon}\right)\left(\sigma_{s}^{1, \epsilon} - \theta_1\right)^2ds
\end{eqnarray}
for some $\sigma_{s,*}^{1, \epsilon}$ lying between $\theta_1$ and $\sigma_{s}^{1, \epsilon}$, with 
\begin{eqnarray*}
\left|\tilde{h}'(\sigma_{s,*}^{1, \epsilon})\right| &\leq& \lambda_1 + \lambda_2\left|\sigma_{s,*}^{1, \epsilon}\right|^m \\
&\leq& \lambda_1 + \lambda_2\left(\left|\sigma_{s}^{1, \epsilon}\right| + \left|\theta_1\right|\right)^m \\
&\leq& \lambda_1 + \lambda_2\left(\left|\sigma_{s}^{1, \epsilon} - \theta_1\right| + 2\left|\theta_1\right|\right)^m\end{eqnarray*}
for some $\lambda_1, \, \lambda_2 > 0$ and some $m \in \mathbb{N}$, which allows us to bound the RHS of (\ref{polybound}) by a linear combination of terms of the form $\Vert\sigma_{\cdot}^{1, \epsilon} - \theta_1\Vert _{L^p(\Omega \times [0, \, T])}^p$ which are all $\mathcal{O}(\epsilon)$ as $\epsilon \to 0^+$ by Proposition~\ref{lem1}. The proof of the theorem is now complete.
\qed
\end{proof}

\begin{proof}[\textbf{Corollary 3.3}]
Let $\mathcal{E}_{t, \, \epsilon} = \{\omega \in \Omega: \, \mathbb{P}[X_t^{1, \epsilon} > 0 \, | \, W^{0}, \, B^{0}, \, \mathcal{G}] > x\}$ for $\epsilon > 0$, $\mathcal{E}_{t, \, 0} = \{\omega \in \Omega: \, \mathbb{P}[X_t^{1, *} > 0 \, | \, W^{0}, \, \mathcal{G}] > x\}$, and observe that
\begin{eqnarray}\label{roc1}
E(x, T) &=& \int_{0}^{T}\big|\mathbb{P}\left[\mathcal{E}_{t, \, \epsilon}\right] - \mathbb{P}\left[\mathcal{E}_{t, \, 0}\right]\big|dt, \nonumber \\
&=& \int_{0}^{T}\Big|\mathbb{P}\big[\mathcal{E}_{t, \, \epsilon} \cap \mathcal{E}_{t, \, 0}^{c}\big] - \mathbb{P}\left[\mathcal{E}_{t, \, 0} \cap \mathcal{E}_{t, \, \epsilon}^{c}\right]\Big|dt, \nonumber \\
&\leq& \int_{0}^{T}\mathbb{P}\left[\mathcal{E}_{t, \, \epsilon} \cap \mathcal{E}_{t, \, 0}^{c}\right]dt + \int_{0}^{T}\mathbb{P}\left[\mathcal{E}_{t, \, 0} \cap \mathcal{E}_{t, \, \epsilon}^{c}\right]dt.  
\end{eqnarray}
Next, for any $\eta > 0$ we have 
\begin{eqnarray}\label{roc2}
&& \mathbb{P}\left[\mathcal{E}_{t, \, \epsilon} \cap \mathcal{E}_{t, \, 0}^{c}\right] \nonumber \\
&& \qquad = \mathbb{P}\bigg[\mathbb{P}\Big[X_t^{1, \epsilon} > 0 \, \Big| \, W^{0}, \, B^{0}, \, \mathcal{G}\Big] > x \geq \mathbb{P}\Big[X_t^{1, *} > 0 \, \Big| \, W^{0}, \, \mathcal{G}\Big]\bigg] \nonumber \\
&& \qquad = \mathbb{P}\bigg[\mathbb{P}\Big[X_t^{1, \epsilon} > 0 \, \Big| \, W^{0}, \, B^{0}, \, \mathcal{G}\Big] > x \nonumber \\
&& \qquad \qquad \qquad > x - \eta > \mathbb{P}\Big[X_t^{1, *} > 0 \, \Big| \, W^{0}, \, \mathcal{G}\Big]\bigg] \nonumber \\
&& \qquad \qquad + \mathbb{P}\bigg[\mathbb{P}\Big[X_t^{1, \epsilon} > 0 \, \Big| \, W^{0}, \, B^{0}, \, \mathcal{G}\Big] > x \nonumber \\
&& \qquad \qquad \qquad \qquad \geq \mathbb{P}\Big[X_t^{1, *} > 0 \, \Big| \, W^{0}, \, \mathcal{G}\Big] \geq x - \eta \bigg] \nonumber \\
&& \qquad \leq \mathbb{P}\bigg[\bigg|\mathbb{P}\Big[X_t^{1, \epsilon} > 0 \, \Big| \, W^{0}, \, B^{0}, \, \mathcal{G}\Big] - \mathbb{P}\Big[X_t^{1, *} > 0 \, \Big| \, W^{0}, \, \mathcal{G}\Big]\bigg| > \eta \bigg] \nonumber \\
&& \qquad \qquad + \mathbb{P}\bigg[x \geq \mathbb{P}\Big[X_t^{1, *} > 0 \, \Big| \, W^{0}, \, \mathcal{G}\Big] \geq x - \eta \bigg] \nonumber \\
&& \qquad \leq \frac{1}{\eta^2}\mathbb{E}\Bigg[\bigg(\mathbb{P}\Big[X_t^{1, \epsilon} > 0 \, \Big| \, W^{0}, \, B^{0}, \, \mathcal{G}\Big] - \mathbb{P}\Big[X_t^{1, *} > 0 \, \Big| \, W^{0}, \, \mathcal{G}\Big]\bigg)^2 \Bigg] \nonumber \\
&& \qquad \qquad + \mathbb{P}\bigg[x \geq \mathbb{P}\Big[X_t^{1, *} \geq 0 \, \Big| \, W^{0}, \, \mathcal{G}\Big] \geq x - \eta \bigg]
\end{eqnarray}
and if we denote by $\mathcal{S}$ the $\sigma$-algebra generated by the volatility paths, since $X_{t}^{1, *}$ is independent of $\mathcal{S}$ and the path of $B^{0}$, by using the Cauchy-Schwarz inequality we find that
\begin{eqnarray*}
&& \int_{0}^{T}\mathbb{E}\Bigg[\bigg(\mathbb{P}\Big[X_t^{1, \epsilon} > 0 \, \Big| \, W^{0}, \, B^{0}, \, \mathcal{G}\Big] - \mathbb{P}\Big[X_t^{1, *} > 0 \, \Big| \, W^{0}, \, \mathcal{G}\Big]\bigg)^2 \Bigg]dt \nonumber \\
&& \qquad = \int_{0}^{T}\mathbb{E}\Bigg[\mathbb{E}\bigg[\mathbb{P}\Big[X_t^{1, \epsilon} > 0 \, \Big| \, W^{0}, \, C_{1}', \, \mathcal{S}, \, \mathcal{G}\Big]  \nonumber \\
&& \qquad \qquad -  \mathbb{P}\Big[X_t^{1, *} > 0 \, \Big| \, W^{0}, \, C_{1}', \, \mathcal{G}\Big] \, \bigg| \, W^{0}, \, B^{0}, \, \mathcal{G} \bigg]^2 \Bigg]dt \nonumber \\
&& \qquad \leq \int_{0}^{T}\mathbb{E}\Bigg[\mathbb{E}\bigg[\bigg(\mathbb{P}\Big[X_t^{1, \epsilon} > 0 \, \Big| \, W^{0}, \, C_{1}', \, \mathcal{S}, \, \mathcal{G}\Big] \nonumber \\
&& \qquad \qquad - \mathbb{P}\Big[X_t^{1, *} > 0 \, \Big| \, W^{0}, \, C_{1}', \, \mathcal{G}\Big]\bigg)^2 \, \bigg| \, W^{0}, \, B^{0}, \, \mathcal{G} \bigg] \Bigg]dt \nonumber \\
&& \qquad = \int_{0}^{T}\mathbb{E}\Bigg[\bigg(\mathbb{P}\Big[X_t^{1, \epsilon} > 0 \, \Big| \, W^{0}, \, C_{1}', \, \mathcal{S}, \, \mathcal{G}\Big] - \mathbb{P}\Big[X_t^{1, *} > 0 \, \Big| \, W^{0}, \, C_{1}', \, \mathcal{G}\Big]\bigg)^2 \Bigg]dt \nonumber \\
&& \qquad = \left\Vert v^{0, \epsilon}\left(\cdot, 0\right) - v^{0}\left(\cdot, 0\right) \right\Vert_{L^2\left(\Omega \times [0, \, T] \right)}^2 \nonumber \\
&& \qquad = \mathcal{O}\left(\epsilon\right),
\end{eqnarray*}
where the last follows by using Morrey's inequality in dimension 1 (see e.g. Evans \cite{EVAN}) and Theorem~\ref{thm5.2}. On the other hand, since $\mathbb{P}[X_t^{1, *} > 0 \, | \, W^{0}, \, \mathcal{G}]$ has a bounded density near $x$, uniformly in $t \in [0, \, T]$, we have
\begin{eqnarray*}
\int_0^T\mathbb{P}\bigg[x \geq \mathbb{P}\Big[X_t^{1, *} > 0 \, \Big| \, W^{0}, \, \mathcal{G}\Big] \geq x - \eta \bigg]dt = \mathcal{O}\left(\eta \right).
\end{eqnarray*}
Therefore, (\ref{roc2}) gives
$\int_0^T\mathbb{P}[\mathcal{E}_{t, \, \epsilon} \cap \mathcal{E}_{t, \, 0}^{c}]dt \leq \frac{1}{\eta^2}\mathcal{O}(\epsilon) + \mathcal{O}(\eta )
$ for any $\eta > 0$, and in a similar way we can obtain
$\int_0^T\mathbb{P}[\mathcal{E}_{t, \, 0} \cap \mathcal{E}_{t, \, \epsilon}^{c}]dt \leq \frac{1}{\eta^2}\mathcal{O}(\epsilon) + \mathcal{O}(\eta )$. Using these two expressions in (\ref{roc1}) and taking $\eta = \epsilon^{p}$ for some $p > 0$, 
we finally obtain
\begin{eqnarray*}
E(x, \, T) \leq \mathcal{O}\left(\epsilon^{p}\right) + \mathcal{O}\left(\epsilon^{1 - 2p}\right),
\end{eqnarray*}
which becomes optimal as $\epsilon \rightarrow 0^{+}$ when $1 - 2p = p \Leftrightarrow p = \frac{1}{3}$. This gives $E(x, \, T) = \mathcal{O}(\epsilon^{\frac{1}{3}})$ as $\epsilon \rightarrow 0^{+}$.
\end{proof}


\appendix
\normalsize
\section{APPENDIX: Proofs of positive reccurence results}

In this Appendix we prove Proposition~\ref{thm3} and Proposition~\ref{thm4}. Both proofs are based on Theorem~2.5 from Bhattacharya and Ramasubramanian \cite{RBSR}, which gives sufficient conditions for an $n$-dimensional Markov process $X$ with infinitesimal generator
\begin{eqnarray*}
L_t = \frac{1}{2}\sum_{i,j = 1}^{n}{a_{i,j}(t, x)\frac{\partial^2}{\partial x_ix_j}} + \sum_{i = 1}^{n}{b_{i}(t, x)\frac{\partial}{\partial x_i}}
\end{eqnarray*}
to be positive reccurent, i.e possess an invariant probability distribution $v$ on $\mathbb{R}^n$ such that
\begin{eqnarray*}
\lim_{T \to +\infty}{\frac{1}{T}\int_{0}^{T}{f\left(X_s\right)ds}} = \int_{\mathbb{R}^n}{f \cdot dv}
\end{eqnarray*}
for any $v$-integrable function $f$. That theorem involves the functions
\begin{eqnarray*}
A_z(s, x) &:=& \sum_{i,j=1}^{n}{a_{i,j}(s, x)\frac{(x_i - z_i)(x_j - z_j)}{\left|x - z\right|_{2}^2}} \\
B(s, x) &:=& \sum_{i=1}^{n}{a_{i,i}(s, x)} \\
C_z(s, x) &:=& 2\sum_{i,j=1}^{n}{b_{i}(s, x)(x_i - z_i)} \\
\underline{\alpha}_z(r) &:=& \inf_{s \geq 0, |x - z|_2 = r}{A_z(s, x)} \\
\overline{\beta}(r; z, t_0) &:=& \sup_{s \geq t_0, |x - z|_2 = r}{\frac{B(s,x) - A_z(s,x) + C_z(s,x)}{A_z(s,x)}} \\
\overline{I}_{z, r_0}(r) &:=& \int_{r_0}^{r}\frac{\overline{\beta}(u;z,0)}{u}du
\end{eqnarray*}
and the conditions implying positive reccurence are the following:
\begin{enumerate}
    \item $a_{i,j}(\cdot, \cdot)$ and $b_i(\cdot, \cdot)$ are Borel measurable on $[0, \, T] \times \mathbb{R}^n$ and bounded on compacts.
    \item For each $N > 0$, there exists a $\delta_N(r) \downarrow 0$ as $r \downarrow 0$ such that for all $t \geq 0$ and $x,y \in \mathbb{R}^n$ with $t, |x|_2, |y|_2 \leq N$ we have 
    \[\left\Vert a_{\cdot, \cdot}(t,x) - a_{\cdot, \cdot}(t,y)\right\Vert_{2} \leq \delta_N(|x - y|_2),\]
    where $\Vert \cdot \Vert_2$ stands for the matrix $2$-norm.
    \item For any compact $K \subset \mathbb{R}^n$ and every $z' \in \mathbb{R}^k$, the function $\frac{B + C_{z'}}{A_{z'}}$ is bounded away from $+\infty$ on $[0, \, +\infty] \times K$.
    \item There exist $z \in \mathbb{R}^n$ and $r_0 > 0$ such that:
    \[
        \int_{r_0}^{+\infty}{e^{-\overline{I}_{z, r_0}(r)}dr} = +\infty
    \]
    and
    \[
        \int_{r_0}^{+\infty}{\frac{1}{\underline{\alpha}_{z}(r)}e^{\overline{I}_{z, r_0}(r)}dr} < +\infty
    \]
\end{enumerate}

We proceed now to our proofs, where we will establish positive reccurence results by showing that the above conditions are satisfied.

\begin{proof}[\bf{Proposition 2.2}] It suffices to show that the two-dimensional continuous Markov process $(\sigma^{1,1}, \, \sigma^{2,1})$ is positive recurrent. To do this, we set $H^i(x) = \int_{0}^{x}\frac{1}{v_ig(y)}dy$ which is a strictly increasing bijection from $\mathbb{R}$ to itself, and then $Z^i = H^i(\sigma^{i,1})$, for $i \in \{1, \, 2\}$. It suffices to show that the two-dimensional process $Z = (Z^{1}, \, Z^{2})$ is positive recurrent. The infinitesimal generator $L_{Z}$ of $Z$ maps any smooth function $F: \mathbb{R}^2 \rightarrow \mathbb{R}$ to
\begin{eqnarray}
L_{Z}F(x, \, y) &=& V^1(x)F_x(x, \, y) + V^2(y)F_y(x, \, y) \nonumber \\
&& + \frac{1}{2}\big(F_{xx}(x, \, y) + F_{yy}(x, \, y)\big) + \lambda F_{xy}(x, \, y) \nonumber
\end{eqnarray}
for $\lambda = \rho_{2,1}\rho_{2,2} < 1$ and $V^i(x) = \frac{\kappa_i(\theta_i - (H^i)^{-1}(x))}{v_ig((H^i)^{-1}(x))} - \frac{v_i}{2}g'((H^i)^{-1}(x))$, with $V^i$ being a continuous and strictly decreasing bijection from $\mathbb{R}$ to itself for $i \in \{1, \, 2\}$. 

We can compute
\begin{eqnarray*}
A_{(z, \, w)}\big(s, \, (x, \, y)\big) &=& \frac{1}{2} + \lambda\frac{(x - z)(y - w)}{(x - z)^2 + (y - w)^2} \nonumber \\
&\geq& \frac{1}{2} + \lambda\frac{-\frac{1}{2}\left((x - z)^2+(y - w)^2\right)}{(x - z)^2 + (y - w)^2} \nonumber \\
&=& \frac{1}{2}(1 - \lambda) > 0
\end{eqnarray*}
and also $B(s, \, (x, \, y)) = 1$ and
\[
C_{(z, \, w)}\big(s, \, (x, \, y)\big) = 2\big(V^1(x)(x - z) + V^2(y)(y - w)\big) 
\]
for all $(x, \, y), (z, \, w) \in \mathbb{R}^2$.
Since the coefficients of $L_{Z}$ are continuous, with the higher order ones being constant, we can easily verify conditions $1.$ and $2$. Moreover, since $B$ and $C_{(z, \, w)}$ are constant in $t$ and continuous in $(x, \, y)$ while $A_{(z, \, w)}$ is lower-bounded by $\frac{1}{2}(1 - \lambda) > 0$, it follows that we have $3.$ as well. 

Next, we choose $z$ and $w$ to be the unique roots of $V^1(x)$ and $V^2(y)$ respectively and we have 
\begin{equation}\label{alpha1}
\underline{\alpha}_{(z, \, w)}(r) = \inf_{(x - z)^2 + (y - w)^2 = r^2}A_{(z, \, w)}\big(s, \, (x, \, y)\big) \geq \frac{1}{2}(1 - \lambda) > 0
\end{equation}
and
\begin{eqnarray}\label{beta}
&&\overline{\beta}\big(r; (z, \, w), 0\big) \nonumber \\
&& \qquad = \sup_{(x - z)^2 + (y - w)^2 = r^2}\frac{B\big(s, \, (x, \, y)\big)-A_{(z, \, w)}\big(s, \, (x, \, y)\big)+C_{(z, \, w)}\big(s, \, (x, \, y)\big)}{A_{(z, \, w)}\big(s, \, (x, \, y)\big)} \nonumber \\
&& \qquad \leq \frac{2}{1 - \lambda} - 1 + \frac{2}{1 + \lambda}\sup_{(x - z)^2 + (y - w)^2 = r^2}{C_{(z, \, w)}\big(s, \, (x, \, y)\big)} \nonumber \\
\end{eqnarray}
since $C_{(z, \, w)}(s, \, (x, \, y))$ is never greater than zero and 
\begin{eqnarray}
A_{(z, \, w)}\big(s, \, (x, \, y)\big) &=& \frac{1}{2} + \lambda\frac{(x - z)(y - w)}{(x - z)^2 + (y - w)^2} \nonumber \\
&\leq& \frac{1}{2} + \lambda\frac{\frac{1}{2}\left((x - z)^2+(y - w)^2\right)}{(x - z)^2 + (y - w)^2} \nonumber \\
&=& \frac{1}{2}(1 + \lambda). \nonumber
\end{eqnarray}
Now fix an $r_0 > 0$ and take any $r > r_0$. For the pair $(x, \, y)$ for which the supremum of $C_{(z, \, w)}(s, \, (x, \, y))$ is attained when $(x - z)^2 + (y - w)^2 = r^2$, we have $x = z + r\cos(\phi_r)$ and $y = w + r\sin(\phi_r)$ for some angle $\phi_r$. Then, we have either $|\cos(\phi_r)| \geq \frac{\sqrt{2}}{2}$ or $|\sin(\phi_r)| \geq \frac{\sqrt{2}}{2}$. If $\cos(\phi_r) \geq \frac{\sqrt{2}}{2}$ holds, we can estimate
\begin{eqnarray}
C_{(z, \, w)}\big(s, \, (x, \, y)\big) &=& 2r\cos(\phi_r)V^1\big(z + r\cos(\phi_r)\big) \nonumber \\
&& \qquad + 2r\sin(\phi_r)V^2\big(w + r\sin(\phi_r)\big) \nonumber \\
&\leq& 2r\cos(\phi_r)V^1\big(z + r\cos(\phi_r)\big) \nonumber \\
&\leq& c_1r \nonumber
\end{eqnarray}
with $c_1 = \sqrt{2}V^1(z + r_0\frac{\sqrt{2}}{2}) < 0$. In a similar way, by using the fact that both $V^1$ and $V^2$ are strictly decreasing, we can find constants $c_2, c_3, c_4 < 0$ such that $C_{(z, \, w)}(s, \, (x, \, y)) < c_2r$, $C_{(z, \, w)}(s, \, (x, \, y)) < c_3r$ and $C_{(z, \, w)}(s, \, (x, \, y)) < c_4r$, when $\cos(\phi_r) \leq -\frac{\sqrt{2}}{2}$, $\sin(\phi_r) \geq \frac{\sqrt{2}}{2}$ and $\sin(\phi_r) \leq -\frac{\sqrt{2}}{2}$ respectively. Thus, for $c^{*} = \max\{c_1, c_2, c_3, c_4\} < 0$ we have $C_{(z, \, w)}(s, \, (x, \, y)) < c_{*}r$, which can be used in (\ref{beta}) to give the estimate
\begin{eqnarray*}
\overline{\beta}\big(r; (z, \, w), 0\big) &\leq& \frac{2}{1 - \lambda} - 1 + \frac{2c^{*}}{1 + \lambda}r
\end{eqnarray*}
for all $r \geq r_0$. This means that for $r_0$ large enough we have
\begin{eqnarray}
\overline{I}_{(z, \, w), r_0}(r) &\leq& \int_{r_0}^{r}\frac{1}{r'}\left(\frac{2}{1 - \lambda} - 1 + \frac{2c^{*}}{1 + \lambda}r'\right)dr' \nonumber \\
&\leq& c^{**}(r - r_0) \nonumber
\end{eqnarray}
for some $c^{**} < 0$ and all $r \geq r_0$. This implies that 
\[
\int_{r_0}^{+\infty}e^{-\overline{I}_{(z, \, w), r_0}(r)}dr = +\infty
\]
and combined with (\ref{alpha1}), it also gives
\[
\int_{r_0}^{+\infty}\frac{1}{\underline{\alpha}_{(z, \, w)}(r)}e^{\overline{I}_{(z, \, w), r_0}(r)}dr \leq \frac{2}{1 - \lambda}\int_{r_0}^{+\infty}e^{c^{**}r}dr < \infty. 
\]
Therefore, we have that all the required conditions are satisfied for the process $Z = (Z^{1}, \, Z^{2})$, which means that $(Z^{1}, \, Z^{2})$ is a positive recurrent diffusion, and thus $(\sigma^{1,1}, \, \sigma^{2,1})$ is positive recurrent as well. 
\qed
\end{proof}

\begin{proof}[\bf{Proposition 2.3}]
We will show first that each volatility process never hits zero. Recalling the standard properties of the scale function $S(x)$ of $\sigma^{1,1}$ (see e.g. Rogers and Williams \cite{RW00}), we have that
\begin{eqnarray*}
S(x) = \int_{\theta_1}^{x}e^{-\int_{\theta_1}^{y}\frac{2\kappa_1(\theta_1 - z)}{v_1^2z\tilde{g}^2(z)}dz}dy
\end{eqnarray*}
and we need to show that $\lim_{n \rightarrow +\infty}S(\frac{1}{n}) = -\infty$. Since $\sup_{x \in \mathbb{R}}\tilde{g}^2(x) \leq 1$ which is strictly less than $\frac{2\kappa_1\theta_1}{v_1^2}$, for $n \geq \frac{1}{\theta_1}$ we have 
\begin{eqnarray}
S\left(\frac{1}{n}\right) &=& -\int_{\frac{1}{n}}^{\theta_1}e^{\int_{y}^{\theta_1}\frac{2\kappa_1(\theta_1 - z)}{v_1^2z\tilde{g}^2(z)}dz}dy \nonumber \\
&\leq& -\int_{\frac{1}{n}}^{\theta_1}e^{\int_{y}^{\theta_1}\frac{(\theta_1 - z)}{\theta_1z}dz}dy \nonumber \\
&\leq& -\int_{\frac{1}{n}}^{\theta_1}e^{\int_{y}^{\theta_1}\frac{1}{z}dz - \int_{y}^{\theta_1}\frac{1}{\theta_1}dz}dy \nonumber \\
&\leq& -\frac{1}{e}\int_{\frac{1}{n}}^{\theta_1}\frac{\theta_1}{y}dy = -\frac{\theta_1}{e}(\ln n + \ln\theta_1) \nonumber
\end{eqnarray}
which tends to $-\infty$ as $n \rightarrow +\infty$. This shows that our volatility processes remain positive forever. 

As our volatility processes are strictly positive, we can set $Z^i = \ln\sigma^{i,1}$ for $i \in \{1, \, 2\}$, and we need to show that $(Z^{1}, \, Z^{2})$ is a positive recurrent diffusion. Again, we can easily determine the infinitesimal generator $L_{Z}$ of $Z = (Z^{1}, \, Z^{2})$, which this time maps any smooth function $F: \mathbb{R}^2 \rightarrow \mathbb{R}$ to
\begin{eqnarray}
L_{Z}F(x, \, y) &=& V^1(x)F_x(x, \, y) + V^2(y)F_y(x, \, y) \nonumber \\
&& + \frac{v_1^2e^{-x}\tilde{g}^2(e^{x})}{2}F_{xx}(x, \, y) + \frac{v_2^2e^{-y}\tilde{g}^2(e^{y})}{2}F_{yy}(x, \, y) \nonumber \\
&& + \lambda v_1 v_2e^{-\frac{x+y}{2}}\tilde{g}(e^{y})\tilde{g}(e^{y})F_{xy}(x, \, y) \nonumber
\end{eqnarray}
for $\lambda = \rho_{2,1}\rho_{2,2} < 1$ and $V^i(x) = e^{-x}(\kappa_i\theta_i - \frac{v_i^2}{2}\tilde{g}^2(e^x)) - \kappa_i$ for $i \in \{1, \, 2\}$, which are again two continuous and strictly decreasing bijections from $\mathbb{R}$ to itself. This can be shown by using the fact that $\tilde{g}$ is increasing and upper-bounded by $1$. 

Using the inequality $ab \geq -\frac{a^2 + b^2}{2}$ we obtain
\begin{eqnarray*}
A_{(z, \, w)}\big(s, \, (x, \, y)\big) &=& \frac{1}{2}\left(\frac{v_1^2e^{-x}\tilde{g}^2(e^{x})(x - z)^2}{(x - z)^2 + (y - w)^2} + \frac{v_2^2e^{-y}\tilde{g}^2(e^{y})(y - w)^2}{(x - z)^2 + (y - w)^2}\right) \nonumber \\
&& + \lambda\frac{v_1v_2e^{-\frac{x+y}{2}}\tilde{g}(e^{x})\tilde{g}(e^{y})(x - z)(y - w)}{(x - z)^2 + (y - w)^2} \nonumber \\
&\geq& \frac{1 - \lambda}{2}\left(\frac{v_1^2e^{-x}\tilde{g}^2(e^{x})(x - z)^2}{(x - z)^2 + (y - w)^2} + \frac{v_2^2e^{-y}\tilde{g}^2(e^{y})(y - w)^2}{(x - z)^2 + (y - w)^2}\right) \nonumber \\
&\geq& \frac{(1 - \lambda)\min\{v_1^2, v_2^2\}}{2}\min\{e^{-x}\tilde{g}^2(e^{x}), \, e^{-y}\tilde{g}^2(e^{y})\} 
\end{eqnarray*}
which is strictly positive. Moreover, we can compute
\begin{eqnarray*}
B\big(s, \, (x, \, y)\big) = \frac{v_1^2e^{-x}\tilde{g}^2(e^{x})}{2} + \frac{v_2^2e^{-y}\tilde{g}^2(e^{y})}{2}
\end{eqnarray*}
and
\[
C_{(z, \, w)}\big(s, \, (x, \, y)\big) = 2\left(V^1(x)(x - z) + V^2(y)(y - w)\right) 
\]
for all $(x, \, y), (z, \, w) \in \mathbb{R}^2$. Since the coefficients of $L_{Z}$ are continuous and $A_{(z, \, w)}(s, \, (x, \, y))$ is strictly positive, it follows that conditions $1.$ and $3.$ are satisfied. To verify $2.$, we pick an $N > 0$ and $x, y, \bar{x}, \bar{y} \in [-N, \, N]$, we set
\begin{eqnarray}
M(x, \, y) &=& \left[\begin{array}{cc}
\frac{v_1^2e^{-x}\tilde{g}^2(e^{x})}{2} & \frac{\lambda v_1 v_2e^{-\frac{x+y}{2}}\tilde{g}(e^{x})\tilde{g}(e^{y})}{2} \nonumber \\
\frac{\lambda v_1 v_2e^{-\frac{x+y}{2}}\tilde{g}(e^{x})\tilde{g}(e^{y})}{2} & \frac{v_2^2e^{-y}\tilde{g}^2(e^{y})}{2}
\end{array} \right] \nonumber 
\end{eqnarray}
and we compute
\begin{eqnarray*}
&& \left\Vert M(x, \, y) - M(\bar{x}, \, \bar{y}) \right\Vert_2 \nonumber \\
&& \qquad = \left(\frac{v_1^2e^{-x}\tilde{g}^2(e^{x})}{2} - \frac{v_1^2e^{-\bar{x}}\tilde{g}^2(e^{\bar{x}})}{2}\right)^2 + \left(\frac{v_2^2e^{-y}\tilde{g}^2(e^{y})}{2} - \frac{v_2^2e^{-\bar{y}}\tilde{g}^2(e^{\bar{y}})}{2}\right)^2 \nonumber \\
&& \qquad \qquad + \left(\frac{\lambda v_1 v_2e^{-\frac{x+y}{2}}\tilde{g}(e^{x})\tilde{g}(e^{y})}{2} - \frac{\lambda v_1 v_2e^{-\frac{\bar{x}+\bar{y}}{2}}\tilde{g}(e^{\bar{x}})\tilde{g}(e^{\bar{y}})}{2}\right)^2 \nonumber \\
&& \qquad \leq C_N\left\| (x, \, y) - (\bar{x}, \, \bar{y}) \right\|_{\mathbb{R}^2}^2,
\end{eqnarray*}
where we have used the two-dimensional mean value theorem on each of the three terms, and the fact that all of the functions involved have a bounded gradient in $[-N, \, N]^2$ (since $\tilde{g}$ has continuous derivatives). Thus, if we take $\delta_N(r) = C_Nr$, we have $2.$ as well.

Next, for some $r_0 > 0$ and all $r \geq r_0$, we compute
\begin{eqnarray}\label{beta3}
&&\overline{\beta}\big(r; (z, \, w), 0\big) \nonumber \\ 
&& \qquad = \sup_{(x - z)^2 + (y - w)^2 = r^2}\frac{B\big(s, \, (x, \, y)\big)-A_{(z, \, w)}\big(s, \, (x, \, y)\big)+C_{(z, \, w)}\big(s, \, (x, \, y)\big)}{A_{(z, \, w)}\big(s, \, (x, \, y)\big)} \nonumber \\
&& \qquad = - 1 + \sup_{(x - z)^2 + (y - w)^2 = r^2}\frac{B\big(s, \, (x, \, y)\big)+C_{(z, \, w)}\big(s, \, (x, \, y)\big)}{A_{(z, \, w)}\big(s, \, (x, \, y)\big)}
\end{eqnarray}
where again we choose $z$ and $w$ to be the unique roots of $V^1(x)$ and $V^2(y)$ respectively. Then, by setting $x = z + r\cos(\phi_r)$ and $y = w + r\sin(\phi_r)$ with $\phi_r \in [0, \, 2\pi]$ for the $(x, \, y)$ for which the above supremum is attained, since $\tilde{g}$ is increasing we have
\begin{eqnarray}\label{est1}
C_{(z, \, w)}\big(s, \, (x, \, y)\big) &=& 2\Bigg(e^{-x}\left(\kappa_1\theta_1 - \frac{v_1^2}{2}\tilde{g}^2(e^x)\right) - \kappa_1\Bigg)(x - z) \nonumber \\
&& + 2\Bigg(e^{-y}\left(\kappa_2\theta_2 - \frac{v_2^2}{2}\tilde{g}^2(e^y)\right) - \kappa_2\Bigg)(y - w) \nonumber \\
&\leq& 2\Bigg(e^{-x}\left(\kappa_1\theta_1 - \frac{v_1^2}{2}\tilde{g}^2(e^z)\right) - \kappa_1\Bigg)(x - z) \nonumber \\
&& + 2\Bigg(e^{-y}\left(\kappa_2\theta_2 - \frac{v_2^2}{2}\tilde{g}^2(e^w)\right) - \kappa_2\Bigg)(y - w) \nonumber \\
&=& 2\kappa_1\left(e^{z-x} - 1\right)(x - z) + 2\kappa_2\left(e^{w-y} - 1\right)(y- w) \nonumber \\
&\leq& 2\min\{\kappa_1, \kappa_2\}\Big(\left(e^{z-x} - 1\right)(x - z) + \left(e^{w-y} - 1\right)(y- w)\Big) \nonumber \\
&=& \kappa r\bigg(\left(e^{-r\cos(\phi_r)} - 1\right)\cos(\phi_r) + \left(e^{-r\sin(\phi_r)} - 1\right)\sin(\phi_r)\bigg) \nonumber \\
\end{eqnarray}
for $\kappa = 2\min\{\kappa_1, \kappa_2\}$, and since $\tilde{g}$ is bounded, for $\xi = \frac{\max\{v_1, \, v_2\}^2}{2}\sup_{x \in \mathbb{R}}\tilde{g}(x)$ we can also show that
\begin{eqnarray}\label{est2}
&& A_{(z, \, w)}\big(s, \, (x, \, y)\big) \nonumber \\
&& \qquad = \xi\left(e^{-x}\frac{(x - z)^2}{(x - z)^2 + (y - w)^2} + e^{-y}\frac{(y - w)^2}{(x - z)^2 + (y - w)^2}\right) \nonumber \\
&& \qquad = \xi\left(e^{-z-r\cos(\phi_r)}\cos^2(\phi_r) + e^{-w-r\sin(\phi_r)}\sin^2(\phi_r)\right) \nonumber \\
&& \qquad = \xi\bigg(e^{-z}\left(e^{-r\cos(\phi_r)}-1\right)\cos^2(\phi_r) + e^{-w}\left(e^{-r\sin(\phi_r)}-1\right)\sin^2(\phi_r)\bigg) \nonumber \\
&& \qquad \qquad + \xi\left(e^{-z}\cos^2(\phi_r) + e^{-w}\sin^2(\phi_r)\right) \nonumber \\
&& \qquad \leq -\xi\bigg(e^{-z}\left(e^{-r\cos(\phi_r)}-1\right)\cos(\phi_r) + e^{-w}\left(e^{-r\sin(\phi_r)}-1\right)\sin(\phi_r)\bigg) \nonumber \\
&& \qquad \qquad + \xi\left(e^{-z} + e^{-w}\right),
\end{eqnarray}
where we have also used the elementary inequality $(e^{ab} - 1)a^2 \leq -(e^{ab} - 1)a$ for $|a| \leq 1$ and $b < 0$. Using (\ref{est1}) and (\ref{est2}) we obtain
\begin{eqnarray}\label{est3}
\frac{C_{(z, \, w)}\big(s, \, (x, \, y)\big)}{A_{(z, \, w)}\big(s, \, (x, \, y)\big)} \leq -r\frac{\kappa}{\xi}\frac{\ell(r)}{\ell(r) + \xi\left(e^{-z} + e^{-w}\right)}
\end{eqnarray}
where
\begin{eqnarray}\label{estlb}
\ell(r) &=& -\xi\bigg(e^{-z}\left(e^{-r\cos(\phi_r)}-1\right)\cos(\phi_r) + e^{-w}\left(e^{-r\sin(\phi_r)}-1\right)\sin(\phi_r)\bigg) \nonumber \\
&\geq& -\xi e^{-z}\left(e^{-r\cos(\phi_r)}-1\right)\cos(\phi_r) \nonumber \\
&=& \xi e^{-z}\left|e^{-r\cos(\phi_r)}-1\right||\cos(\phi_r)| \nonumber \\
&\geq& \xi \min\{e^{-z}, e^{-w}\}\frac{\sqrt{2}}{2}\min\bigg\{\left|e^{-r_0\frac{\sqrt{2}}{2}}-1\right|, \, \left|e^{r_0\frac{\sqrt{2}}{2}}-1\right|\bigg\}  
\end{eqnarray}
since we take $r \geq r_0$, and without loss of generality we can assume that $|\cos(\phi_r)| \geq \frac{\sqrt{2}}{2}$. Thus, (\ref{est3}) implies that there is a universal $c^{*} < 0$ such that
\begin{eqnarray*}
\frac{C_{(z, \, w)}\big(s, \, (x, \, y)\big)}{A_{(z, \, w)}\big(s, \, (x, \, y)\big)} \leq c^{*}r
\end{eqnarray*}
when $r \geq r_0$. Plugging the last in (\ref{beta3}) we obtain
\begin{eqnarray}\label{beta4}
&&\overline{\beta}\big(r; (z, \, w), 0\big) \nonumber \\
&& \qquad = - 1 + \sup_{(x - z)^2 + (y - w)^2 = r^2}\frac{B\big(s, \, (x, \, y)\big)+C_{(z, \, w)}\big(s, \, (x, \, y)\big)}{A_{(z, \, w)}\big(s, \, (x, \, y)\big)} \nonumber \\
&& \qquad \leq - 1 + pc^{*}r + \sup_{(x - z)^2 + (y - w)^2 = r^2}\frac{B\big(s, \, (x, \, y)\big)+(1-p)C_{(z, \, w)}\big(s, \, (x, \, y)\big)}{A_{(z, \, w)}\big(s, \, (x, \, y)\big)} \nonumber \\
\end{eqnarray}
for all $r \geq r_0$ and a $p \in [0, \, 1]$ which will be chosen later. 

We will show now that the last term in the RHS of (\ref{beta4}) above is negative for $r \geq r_0$ with $r_0$ large enough (depending on $p$). Indeed, by using (\ref{est1}), the definition of $B(s, \, (x, \, y))$, and the fact that $\tilde{g}$ is upper-bounded, we can obtain the estimate
\begin{eqnarray}\label{est5}
&&\sup_{(x - z)^2 + (y - w)^2 = r^2}\frac{B\big(s, \, (x, \, y)\big)+(1-p)C_{(z, \, w)}\big(s, \, (x, \, y)\big)}{A_{(z, \, w)}\big(s, \, (x, \, y)\big)} \nonumber \\
&& \qquad \leq \sup_{(x - z)^2 + (y - w)^2 = r^2} \Bigg\{ \frac{\kappa^{*}\big(\left(e^{z - x} - 1\right)(x - z) + \left(e^{w - y} - 1\right)(y - w)\big)}{A_{(z, \, w)}\big(s, \, (x, \, y)\big)} \nonumber \\
&& \qquad \qquad + \frac{\xi(e^{-x}+e^{-y})}{A_{(z, \, w)}\big(s, \, (x, \, y)\big)} \Bigg\}
\end{eqnarray}
where as before, we have $\xi = \frac{\max\{v_1, \, v_2\}^2}{2}\sup_{x \in \mathbb{R}}\tilde{g}(x)$, and $\kappa^{*} = (1 - p)\kappa$. The numerator in the last supremum can easily be shown to tend to $-\infty$ when $x$ or $y$ tends to $\pm \infty$, which happens when $r \rightarrow +\infty$. Thus, for $r \geq r_0$ with $r_0$ large enough, the RHS of (\ref{est5}) is negative. 

The last can be used in (\ref{beta4}) to give
\begin{eqnarray}\label{beta5}
\overline{\beta}\big(r; (z, \, w), 0\big) &\leq& - 1 + pc^{*}r
\end{eqnarray}
for all $r \geq r_0$, with $c^{*} < 0$. On the other hand, we can also compute
\begin{eqnarray}\label{alpha2}
\underline{\alpha}_{(z, \, w)}(r) &=& \inf_{(x - z)^2 + (y - w)^2 = r^2}A_{(z, \, w)}\big(s, \, (x, \, y)\big) \nonumber \\
&\geq& \frac{(1 - \lambda)\min\{v_1^2, v_2^2\}}{2}e^{-\max\{z, \, w\} - r}\tilde{g}^2\left(e^{\max\{z, \, w\} + r}\right) \nonumber \\
\end{eqnarray}
It follows from (\ref{beta5}) that for $r_0$ large enough we have
\begin{eqnarray}\label{est6}
\overline{I}_{(z, \, w), r_0}(r) &\leq& \int_{r_0}^{r}\frac{1}{r'}\left(- 1 + pc^{*}r'\right)dr' \leq pc^{*}(r - r_0)
\end{eqnarray}
with $c^{*} < 0$, for all $r \geq r_0$, which implies that 
\[
\int_{r_0}^{+\infty}e^{-\overline{I}_{(z, \, w), r_0}(r)}dr = +\infty, 
\]
while we can use (\ref{est6}), (\ref{alpha2}) and the fact that $\tilde{g}$ is lower-bounded by something positive to find that
\[
\int_{r_0}^{+\infty}\frac{1}{\underline{\alpha}\big(r; (z, \, w), 0\big)}e^{\overline{I}_{(z, \, w), r_0}(r)}dr < C_{\alpha}\int_{r_0}^{+\infty}e^{(1+pc^{*})r}dr. 
\]
for some $C_{\alpha} > 0$, and it remains to show that the RHS of the above is finite. The last is the case when $pc^{*} < -1$, which is achieved by taking $\eta$ small enough and $r_0$ big enough. Indeed, this forces $e^{-w}$ and $e^{-z}$ to be arbitrarily close to each other, and the lower bound of $\ell$ (given by (\ref{estlb})) to be arbitrarily close to $\frac{\sqrt{2}\xi}{2}e^{-z}$. This brings $c^{*}$ (obtained in (\ref{est3})) arbitrarily close to 
\begin{eqnarray}
-\frac{\kappa}{\xi}\frac{\sqrt{2}}{\sqrt{2}+ 4} = -\frac{\min\{\kappa_1, \kappa_2\}}{\max\{v_1^2, v_2^2\}\sup_{x \in \mathbb{R}}\tilde{g}(x)}\frac{4\sqrt{2}}{\sqrt{2}+ 4} \nonumber 
\end{eqnarray} 
which is less by $-1$ by our initial assumptions. Therefore, if $p$ is chosen to be sufficiently close to $1$, we have $pc^{*} < -1$ as well and this completes the proof.
\end{proof}



\begin{thebibliography}{10}


\bibitem{Aldous} Aldous, D. Exchangeability and related topics (1985), Ecole d\textquoteright Ete St Flour 1983, Springer Lecture Notes in Mathematics, 1117, 1--198

\bibitem{RBSR} Bhattacharya, R. and Ramasubramanian, S. Recurrence and ergodicity of diffusions, \emph{J. Multivariate Anal.}, {\bf 12(1)} (1982), 95--122.

\bibitem{Brezis} Br\'{e}zis, H. \emph{Functional analysis, Sobolev spaces and partial differential equations}. Universitext. New York ; London: Springer (2011)

\bibitem{BR12} Bujok, K. and Reisinger, C. 
Numerical valuation of basket credit derivatives in structural jump-diffusion models, 
\emph{J. Comput. Finance}, {\bf 15} (2012), 115--158.

\bibitem{BHHJR}
Bush, N.; Hambly, B.M.; Haworth, H.; Jin, L., and Reisinger, C. Stochastic evolution equations in portfolio credit modelling, \emph{SIAM J. Financial Math.}, {\bf 2} (2011), 627--664.

\bibitem{CMZ} Cvitanic, J.; Ma, J. and Zhang, J. Law of large numbers for self-exciting correlated defaults, \emph{Stochastic Process. Appl.}, {\bf 122(8)} (2012), 2781--2810.


\bibitem{PRST} Dai Pra, P.; Runggaldier, W.; Sartori, E. and Tolotti, M. Large portfolio losses: A dynamic contagion mode, \emph{Ann. Appl. Probab.}, {\bf 19(1)} (2009), 347--394.

\bibitem{PT09} Dai Pra, P. and Tolotti, M. Heterogeneous credit portfolios and the dynamics of the aggregate losses, \emph{Stochastic Process. Appl.}, {\bf 119(9)} (2009), 2913--2944.

\bibitem{EK05} Ethier, S. and Kurtz, T. \emph{Markov processes: Characterization and convergence}. Wiley series in probability and statistics. Hoboken, N.J.: Wiley-Interscience (2005)  

\bibitem{EVAN} Evans, L. \emph{Partial differential equations}. Graduate studies in mathematics ; v.19. Providence, Rhode Island: American Mathematical Society (2015)  

\bibitem{FOU} Fouque, J.; Papanicolaou, G. and Sircar, K. \emph{Derivatives in financial markets with stochastic volatility}. Cambridge: Cambridge University Press (2000)

\bibitem{GSS2} Giesecke, K; Spiliopoulos, K. and Sowers, R. B. Default clustering in large portfolios: Typical events, \emph{Ann. Appl. Probab.} {\bf 23(1)} (2013), 348--385.

\bibitem{GSS} Giesecke, K; Spiliopoulos, K.; Sowers, R. B. and Sirignano, J. A. Large portfolio
asymptotics for loss from default, \emph{Math. Finance} {\bf 25} (2015), 77--114.

\bibitem{GR} Giles, M.B. and Reisinger, C. Stochastic finite differences and multilevel Monte Carlo for a class of SPDEs in finance, \emph{SIAM J. Financial. Math.}, {\bf 3} (2012), 572--592.

\bibitem{HK17} Hambly, B. and Kolliopoulos, N. Stochastic evolution equations for large portfolios of stochastic volatility models, \emph{SIAM J. Financial Math.} {\bf 8} (2017), 962--1014.

\bibitem{HKE17} Hambly, B. and Kolliopoulos, N. Erratum: Stochastic evolution equations for large portfolios of stochastic volatility models, \emph{SIAM J. Financial Math.} {\bf 10} (2019), 857--876. 

\bibitem{HK19} Hambly, B and Kolliopoulos, N. Stochastic PDEs for large portfolios with general mean-reverting volatility processes. Submitted for publication (2020). Arxiv ID: 1906.05898

\bibitem{HL16} Hambly, B. and Ledger, S. A stochastic McKean--Vlasov equation
for absorbing diffusions on the half-line, \emph{Ann. Appl. Probab.} {\bf 27} (2017), 2698--2752.



\bibitem{HV15}  Hambly, B. and Vaicenavicius, J. The 3/2 model as a stochastic volatility approximation for a large-basket price-weighted index, \emph{Int. J. Theor. Appl. Finance}, {\bf 18} (2015), no. 6, 1550041.




\bibitem{KKP} Kotelenez, P. and Kurtz, T. Macroscopic limits for stochastic partial differential equations of McKean-Vlasov type, \emph{Probab. Theory Related Fields},  {\bf 146} (2010), 189--222

\bibitem{Krylov} Krylov, N. A $W_2^n$ -theory of the Dirichlet problem for SPDEs in general smooth domains, \emph{Probab. Theory Related Fields}, {\bf 98} (1994), 389--421.

\bibitem{KR81} Krylov, N. and Rozovskii, B. Stochastic evolution equations, \emph{J. Soviet Math.} {\bf 16} (1981), 1233--1277.

\bibitem{KX99} Kurtz, T.G. and Xiong, J. Particle representations
For a class of non-linear SPDEs, \emph{Stochastic Process. Appl.}, {\bf 83} (1999), 103--126.

\bibitem{Ledger14}
Ledger, S. Sharp regularity near an absorbing boundary
for solutions to second order SPDEs in a half-line with constant coefficients,
\emph{Stoch. Partial Differ. Equ. Anal. Comput.},
{\bf 2} (2014), 1--26.




\bibitem{PLEV} Levy, P. \emph{horie de l’addition des variables alatoires. 2. ed}. Collection des monographies des probabilits. Paris: Gauthier-Villars (1954)

\bibitem{RW00} Rogers, L.C.G. and Williams, D. \emph{Diffusions, Markov Processes and Martingales}. Cambridge: Cambridge University Press ; 2nd ed. [reissued] (2000)


\bibitem{SG} Sirignano, J.A. and Giesecke, K., Risk Analysis for Large Pools of Loans, 
\emph{Management Science} {\bf 65} (2018), 107--121.

\bibitem{STG} Sirignano, J. A.; Tsoukalas, G.; Giesecke, K. Large-scale loan portfolio selection, \emph{Oper. Res.}, {\bf 64} 
(2016), 1239--1255.

\bibitem{SSG} Spiliopoulos, K.; Sirignano, J. A. and Giesecke, K. 
Fluctuation analysis for the loss from default, \emph{Stochastic Process. Appl.} {\bf 124} (2014), 2322--2362.

\bibitem{SSO1} Spiliopoulos, K.; Sowers, R. B. Recovery rates in investment-grade pools of credit assets: A large deviations analysis, \emph{Stochastic Process. Appl.} {\bf 121(12)} (2011), 2861--2898.

\bibitem{SSO2} Spiliopoulos, K.; Sowers, R. B. Default Clustering in Large Pools: Large Deviations, \emph{SIAM J. Financial Math.} {\bf 6(1)} (2015), 86--116.

\end{thebibliography}

%
%

\end{document}